\documentclass[review,3p,10pt]{elsarticle}
\biboptions{sort&compress}

\usepackage{tabularx}
\usepackage{graphicx}
\usepackage{epstopdf}
\usepackage{amsmath}
\usepackage{cases}
\usepackage[caption=false,font=footnotesize]{subfig}
\usepackage{amssymb}
\usepackage{mathtools}
\usepackage{bm}
\usepackage{multirow}
\usepackage{algorithm}
\usepackage{algorithmic}
\usepackage{amsthm}
\usepackage{mathrsfs}
\usepackage{textcomp,mathcomp}
\usepackage{threeparttable}
\usepackage{bbding}
\usepackage{pifont}
\usepackage{multicol}
\usepackage{enumitem}
\usepackage{adjustbox}
\usepackage{xcolor}

\usepackage{makecell}

\usepackage{lineno}

\usepackage{url}

\newcommand{\tabincell}[2]{\renewcommand\arraystretch{0.8}\begin{tabular}{@{}#1@{}}#2\end{tabular}}

\journal{International Journal of Hydrogen Energy}

\begin{document}
\begin{frontmatter}

\title{Plant-Wide Hierarchical Electricity-Heat 
Coordination for Large-Scale Cold-Region 
ReP2H Plants via Bidirectional Thermal Coupling}

\author[label1]{Yiwei~Qiu}
\author[label1]{Baiping~Zhu}
\author[label1]{Tao~Wu\corref{cor1}}
\ead{taowu.cq@163.com}
\author[label1]{Shi~Chen}
\author[label1]{Buxiang~Zhou}
\author[label1]{Kaigui~Xie}
\cortext[cor1]{Corresponding author}

\address[label1]{College of Electrical Engineering, Sichuan University, Chengdu, 610065, China}

\begin{abstract}
Large-scale renewable power-to-hydrogen (ReP2H) plants in cold regions suffer from prolonged startup and repeated thermal stress during frequent startup-shutdown operation. 
The situation becomes worse due to the lack of coordinated heat management among the alkaline electrolysis stacks, balance of plant (BoP), plant thermal utility system (PTUS), and plant building. This paper presents a plant-wide thermal topology and a hierarchical electricity-heat  management framework to address the issues. Bidirectional thermal coupling between the stack cluster and PTUS enables preheating, thermal standby, and waste heat recovery, while minute-scale production scheduling is coordinated with second-scale thermal regulation. Case studies based on an $80~\mathrm{MW}$ plant in Northern China show that the proposed framework eliminates cold startups in year round, increases hydrogen yield by $1.50\%$, improves energy and exergy efficiencies by $0.99$ and $4.33$ percentage points, respectively, and reduces the levelized cost of hydrogen by $3.22\%$. It also reduces thermal fatigue damage and startup-shutdown-induced voltage degradation.
\end{abstract}


\end{frontmatter}

\section*{Nomenclature}

\addcontentsline{toc}{section}{Nomenclature}
\begin{multicols}{2}

	\footnotesize
	\setlength{\columnsep}{18pt}

\noindent\textbf{Abbreviations}

\begin{description}[leftmargin=2.3cm, labelwidth=2.1cm, labelsep=0.2cm, align=left, itemsep=0pt, topsep=0pt, font=\normalfont]
    \item[AWE] Alkaline water electrolysis
    \item[BHE] Bidirectional heat exchanger
    \item[BoP] Balance of plant
    \item[CWHE] Cooling water heat exchanger
    \item[PTUS] Plant thermal utility system
    \item[LCOH] Levelized cost of hydrogen
    \item[MILP] Mixed-integer linear programming
    \item[ReP2H] Renewable power-to-hydrogen
    \item[SEP] Gas-lye separator
\end{description}

\vspace{5pt}

\noindent\textbf{Indices}

\begin{description}[leftmargin=2.3cm, labelwidth=2.1cm, labelsep=0.2cm, align=left, itemsep=0pt, topsep=0pt, font=\normalfont]
    \item[$i,j$] Stack and group indices
    \item[$k$] Scheduling interval index
    \item[{$q$}] {Lower-layer control-step index}
    \item[$\ell$] Temperature cycle index
    \item[{$\xi$}] {Seasonal scenario index}
    \item[{$y$}] {Year index}
\end{description}

\vspace{5pt}

\noindent\textbf{Variables}

\begin{description}[leftmargin=2.3cm, labelwidth=2.1cm, labelsep=0.2cm, align=left, itemsep=0pt, topsep=0pt, font=\normalfont]
    \item[$b_{i,j,k}^{\mathrm{P/S/I}}$] Production, thermal standby, and idle state indicators for stack $i$ in group $j$
    \item[$b_{i,j,k}^{\mathrm{SU/SD/SP}}$] Indicators of startup, shutdown, and the transition from standby to production
    \item[$b_{\mathrm{hc}}$] BHE heat exchange direction
\item[$I_{i,j,q}$] Stack current
    \item[$I_{i,j,q}^{\mathrm{des}}$]
Desired stack current

\item[$I_{i,j,q}^{\max}$]
Maximum allowable stack current
    \item[{$P_{i,j,k}^{\mathrm{ele}}$}] {Scheduled stack electrolytic power}
 \item[$P_{\mathrm{SB}}$]
Stack standby power

\item[$P_{\text{min/max}}$]
Stack power limits
    \item[$P_k^{\mathrm{RE}}$] Available renewable power
\item[{$P_{k,\mathrm{pump/heat}}^{(\xi)}$}] {Total power of the PTUS and lye circulation pumps and electric boiler input power}
    \item[$P_{\mathrm{boiler}}^{\mathrm{in/out}}$] Electric boiler input and output power

\item[$P_{\mathrm{sys}}$]
Net system power demand
    \item[$F_{i,j,k}$] Linearized hydrogen production rate
\item[$\dot{n}_{\mathrm{H}_2}$] Total hydrogen production rate
\item[$V_{\mathrm{m,N}}$] Molar volume at normal conditions
    \item[$Q_{i,j}^{\mathrm{ele}}$, $Q_{i,j,\mathrm{s,diss}}$] Stack heat generation and dissipation rates
    \item[$Q_{\mathrm{atm/vent}}$] Building envelope and ventilation heat loss
    \item[$T_{i,j,\mathrm{s,out}}$] Stack outlet temperature
    \item[$T_{j,\mathrm{s,in}}$]
Stack inlet temperature

\item[$T_{i,j,k-1}^{\mathrm{actual}}$] Measured feedback stack temperature
    \item[$T_{j,\mathrm{sep,in/out}}$] SEP inlet/outlet temperature
    \item[$T_{j,\mathrm{h,out}}$] BHE outlet temperature
    \item[$T_{j}^{\mathrm{s/r}}$] PTUS supply/return temperatures
    \item[$T_{\mathrm{plant}}$]
Plant indoor air temperature

\item[$T_{\mathrm{atm}}$]
Ambient temperature

\item[$T_{\mathrm{rad/w}}$]
Radiator and water temperatures
\item[$t_{\mathrm{su}}$] Startup duration
    \item[$U_{i,j}^{\mathrm{cell}}$] Cell voltage
    \item[$v_{i,j,\mathrm{lye}}$] Stack lye  flow rate
    \item[$v_{j}^{\mathrm{HS/HD}}$] PTUS source/demand branch flow rate
    \item[$v_{\mathrm{pump}}$] Total PTUS circulation flow rate
    \item[$v_{jk}$] PTUS flow rate from node $j$ to node $k$
    \item[$v_{\mathrm{h}}$] Radiator water flow rate

 \item[$v_{j,\mathrm{hs/hr}}$]
PTUS supply and heat recovery flow rates

\item[$v_{j,\mathrm{pre/rec}}$]
Preheating and recovery branch flow rates
    \item[$v_{j,\mathrm{c}}$] Cooling water flow rate
    \item[$\alpha_{\mathrm{vent}}$] Ventilation opening coefficient
    \item[$\bar{D}$] Average cumulative mechanical damage
    \item[$\Delta U_{i,j,y}$] Annual startup-shutdown-induced voltage increment of stack $i$ in group $j$
    \item[$N_{\mathrm{CS/HS},i,j,y}$] Annual counts of cold/hot startups
    \item[$N_{\mathrm{op},k}$] Number of stacks in production
    \item[$\gamma_{i,j,y}$] Hydrogen conversion rate
    \item[$\dot{E}_{x,m}^{\mathrm{ph/ch}}$] Physical/chemical exergy rates of material stream $m$
  \item[$\dot{E}_{x,\mathrm{AWE}}$] Electrical exergy input to the AWE system
\item[$\dot{E}_{x,\mathrm{boiler}}$] Electrical exergy input to the boiler
\item[$\dot{E}_{x,\mathrm{pump}}$] Electrical exergy input to the pumps
    \item[$\dot{E}_{x,\mathrm{d/product}}$] Exergy destruction/product exergy rates
    \item[$\dot{E}_{x,\mathrm{in,total}}$] Total input exergy rate
    \item[$\eta_{\mathrm{sys/ex}}$]
Instantaneous energy/exergy efficiencies

\item[$\bar{\eta}_{\mathrm{sys/ex}}$]
Average energy/exergy efficiencies
\end{description}

\vspace{5pt}

\noindent\textbf{Parameters}

\begin{description}[leftmargin=2.3cm, labelwidth=2.1cm, labelsep=0.2cm, align=left, itemsep=0pt, topsep=0pt, font=\normalfont]
   \item[$A_{\mathrm{b/rad/water}}$]
    Building envelope, radiator, and water-side heat transfer areas

\item[$A_{\mathrm{c/h}}$]
CWHE and BHE heat transfer areas

\item[$A_{\mathrm{s/sep,diss}}$]
Stack and SEP heat dissipation areas
 \item[{$\alpha_{\mathrm{heat/pump}}^{(\xi)}$}]
 {Seasonal electric boiler and pump power fitting coefficients}

\item[{$b_{\mathrm{heat/pump}}^{(\xi)}$}]
 {Seasonal electric boiler and pump power fitting intercepts}

\item[$C_{i,j,\mathrm{s}}$]
Thermal capacity of stack $i$ in group $j$

\item[$C_{j,\mathrm{sep/he/h/ce/c}}$]
Thermal capacities of the SEP, BHE structure and lye, water in BHE, and lye and cooling water in CWHE in group $j$

\item[$C_{\mathrm{plant}}$]
Thermal capacity of the indoor air and effective building thermal mass

\item[$C_{\mathrm{rad/w}}$]
Thermal capacities of the radiator and the water in it
    \item[$c_{\mathrm{air/h/c/lye}}$] Specific heat capacities of air, PTUS water, cooling water, and lye
    \item[$d,l$] PTUS pipe diameter and length
    \item[$f_1,f_2$, $F$] Faradaic efficiency parameters and Faraday constant
    \item[$H_{i,j,y}$] Annual hydrogen production without considering degradation
    \item[$h_{\mathrm{s/sep}}$]
Stack and SEP heat transfer coefficients

\item[$k_{\mathrm{c/h}}$]
CWHE and BHE heat transfer coefficients

\item[$U_{\mathrm{b/rad/water}}$]
Heat transfer coefficients of the building envelope, radiator, and radiator water side

\item[$U_{\mathrm{vent,max}}$]
Maximum ventilation heat loss coefficient
    \item[$\mathit{HHV}_{\mathrm{H}_2}$] Higher heating value of hydrogen
 \item[$N^{\mathrm{cell}}$]
Number of cells per stack

\item[$N$, $N_{\mathrm{g}}$]
Numbers of stacks per group and groups

\item[$N_{\mathrm{rad}}$]
Number of radiators

\item[$N_{\mathrm{T}}$]
Number of scheduling intervals
    \item[$R_{\mathrm{CWHE/BHE}}$] CWHE and BHE thermal resistances
\item[$T_{\mathrm{in,set}}$]
Group inlet temperature setpoint

\item[{$T_{\mathrm{set}}$}]
{Stack target temperature}
\item[{$T_{j}^{\mathrm{HS,s/r}}$}] {PTUS heat-source supply/return temperatures}
  \item[$t_{\mathrm{preheat}}$]
Preheating time
\item[{$t_{\mathrm{a}}$}] {Thermal management advance time}

\item[$t_{\mathrm{start}}$]
Scheduled startup time
\item[{$\Delta T_{\mathrm{d}}$}] {Temperature deadband of thermal standby}

\item[$T_{\mathrm{upper/lower}}$]
Upper-layer scheduling horizon and lower-layer control period

\item[$\Delta t_{\mathrm{upper/lower}}$]
Upper-layer scheduling step and lower-layer control step
    \item[$U^{\mathrm{rev/th/max}}$] Reversible, thermoneutral, and maximum cell voltages
    \item[$\alpha$, $\beta$] Fatigue model parameters

    \item[$\Delta T_{i,j,\ell}$] Temperature cycle range
    \item[$\Delta U_{\mathrm{c/h}}$] Voltage degradation per cold and hot startup
    \item[$\varepsilon_{\mathrm{s/rad/sep}}$]
Stack, radiator, and SEP emissivities
    \item[{$\varphi_{\mathrm{s}}$}] {Stack diameter}
    \item[$\sigma$] Stefan-Boltzmann constant
    \item[$p$] System pressure
\item[$\eta_{i,j}^{\mathrm{cell}}$]
Faradaic efficiency

\item[$\eta_{\mathrm{boiler/pump}}$]
Electric boiler and PTUS circulation pump efficiencies
    \item[$\psi_{jk}^{\mathrm{HL}}$] Temperature retention factor for pipe $jk$
   \item[$\rho_{\mathrm{air/h/c/lye}}$] Air, PTUS water, cooling water, and lye densities

    \item[$r_1,r_2,r_3,s,$] Semi-empirical cell voltage parameters

    \item[$t_1,t_2,t_3$]
    \item[{$I_{\mathrm{HL}}$}] {High-load current threshold}

\end{description}

\end{multicols}

\section{Introduction}
\label{sec:intro}


Renewable hydrogen is increasingly recognized as an important option for decarbonizing hard-to-abate industrial sectors and linking large-scale renewable electricity with hydrogen-based chemical production \cite{Li2025Redesigning,yang2022breaking}. Global deployment of renewable-powered electrolysis is accelerating as clean-hydrogen projects expand toward industrial scale \cite{IEA2025}. Renewable power-to-hydrogen (ReP2H) therefore provides an important pathway for large-scale renewable energy utilization and low-carbon hydrogen production. However, the efficiency, flexibility, and durability of large-scale ReP2H plants remain major challenges under fluctuating renewable power.

Alkaline water electrolysis (AWE) is widely used in ReP2H projects because of its technological maturity, durability, and relatively low cost \cite{huang2025review}. As green hydrogen projects expand from megawatt to gigawatt scale, renewable power fluctuations impose stricter requirements on operating efficiency, dynamic response, and equipment reliability \cite{dowling2020role}. A severe stack failure requiring major overhaul can increase the levelized cost of hydrogen (LCOH) by about 1\%--3\% \cite{lin2026reliability}. Efficient, flexible, and durable operation is therefore essential for economical large-scale hydrogen production \cite{xiao2020optimal, matute2021multi}.

Thermal management is particularly important for AWE efficiency and reliability \cite{qi2023thermal}. Stack temperature directly affects electrochemical performance and lifetime \cite{ali2016developing, david2019advances}, and AWE systems generally operate efficiently at 343--363 K \cite{kojima2018development}. Many large ReP2H projects are located in renewable-rich cold regions, including Northern China \cite{fan2025economic} and Northern Europe \cite{de2024worldwide, isooja2025life}, where winter temperatures can fall to $-30~^\circ\mathrm{C}$ or below \cite{government2025shenneng, zhang2026alkaline, zhai2024review, meng2024advantages}. Low ambient temperature prolongs cold-startup duration, while repeated power variations cause stack temperature cycling. The resulting thermal stress can accelerate seal damage, material creep, and electrode degradation \cite{dutton2000experience, todd2014thermodynamics, brauns2022experimental}.

Plant scale further complicates thermal management. Electrolysis stacks are commonly arranged in multi-stack shared-BoP configurations. A set of $N$ parallel stacks sharing a common BoP is termed an \emph{$N$-in-1 group} (hereafter referred to as a \emph{group}) \cite{qiu2026dynamic}. Within each group, stacks, gas-lye separators (SEPs), cooling water heat exchangers (CWHEs), and lye circulation loops interact through coupled mass and heat flows. Across the plant, these groups further interact with the plant thermal utility system (PTUS) (i.e., the plant-internal hot-water loop comprising the electric boiler, circulation pump, supply/return pipes, and heating branches) and building thermal loads. Their electrical and thermal dynamics span different time scales. Production scheduling determines stack commitment and power allocation over minutes, whereas temperature and flow regulation must respond within seconds. A single centralized control layer therefore has difficulty achieving both plant-wide economic scheduling and fast thermal regulation.

Existing studies have investigated AWE thermal dynamics, multi-stack scheduling, and electricity-heat integration, as reviewed in Section~\ref{sec:lit_review}. However, most focus on individual electrolyzers, simplified stack clusters, or electrolysis systems coupled to external heating networks. A plant-level framework that coordinates heat generation, recovery, transfer, and demand among the stack cluster, BoP, PTUS, and plant building remains lacking.

\begin{table}[tb]\footnotesize
	\renewcommand{\arraystretch}{1.05}
	\setlength{\tabcolsep}{3pt}
	\caption{Summary of recent studies on energy management for hydrogen production systems.}
	\vspace{6pt}
	\label{tab:literature_review}
	\centering
	\begin{adjustbox}{max width=\textwidth}
		\begin{tabular}{@{}ccccccccc@{}}
			\hline\hline
			\multirow{2}{*}{Literature} & \multirow{2}{*}{\tabincell{c}{Configuration}} & \multicolumn{6}{c}{Considered processes} & \multirow{2}{*}{Method} \\
			\cline{3-8}
			& & \tabincell{c}{Heat transfer\\among stacks} & \tabincell{c}{BoP\\heat transfer} & \tabincell{c}{Stack cluster\\control} & \tabincell{c}{Preheating} & \tabincell{c}{Waste heat\\recovery} & \tabincell{c}{Thermal\\integration} & \\
			\hline
			Jin 2025 \cite{jin2025alkaline} & 1-in-1 & $\times$ & $\checkmark$ & $\times$ & $\times$ & $\times$ & $\times$ & \tabincell{c}{Coupled thermal and\\electrochemical modeling} \\
			Meng 2026 \cite{meng2026model} & 1-in-1 & $\times$ & $\checkmark$ & $\times$ & $\times$ & $\times$ & $\times$ & Model optimization \\
			Zhong 2025 \cite{zhong2025improving} & 1-in-1 & $\times$ & $\checkmark$ & $\times$ & $\checkmark$ & $\checkmark$ & Heat storage tank & \tabincell{c}{Unsteady-state\\thermodynamic modeling} \\
			Guan 2025 \cite{guan2025dynamic} & Multiple 1-in-1 & $\times$ & $\times$ & $\checkmark$ & $\times$ & $\times$ & $\times$ & Rolling optimization \\
			Guan 2026 \cite{guan2026region} & $N$-in-1 & $\times$ & $\checkmark$ & $\checkmark$ & $\times$ & $\times$ & $\times$ & Hierarchical scheduling \\
			Zou 2025 \cite{zou2025control} & $N$-in-1 & $\times$ & $\times$ & $\checkmark$ & $\times$ & $\times$ & $\times$ & Multiple-stack control \\
			Wang 2025 \cite{wang2025collaborative} & Multiple 1-in-1 & $\times$ & $\times$ & $\checkmark$ & $\times$ & $\times$ & $\times$ & Rolling optimization \\
			Firdous 2026 \cite{firdous2025utility} & Multiple 1-in-1 & $\times$ & $\times$ & $\times$ & $\times$ & $\times$ & $\times$ & \tabincell{c}{Multiphysics operational\\modeling} \\
			Chen 2026 \cite{chen2025mean} & Multiple 1-in-1 & $\times$ & $\times$ & $\checkmark$ & $\times$ & $\times$ & $\times$ & Mean-field control \\
			Xu 2025 \cite{xu2025optimization} & Multiple 1-in-1 & $\times$ & $\checkmark$ & $\checkmark$ & $\checkmark$ & $\times$ & $\times$ & Optimization control \\
			Ma 2025 \cite{ma2025cold} & $N$-in-1 & $\times$ & $\times$ & $\checkmark$ & $\checkmark$ & $\times$ & $\times$ & Startup and shutdown control \\
			Qiu 2023 \cite{qiu2023extended} & Multiple 1-in-1 & $\times$ & $\times$ & $\checkmark$ & $\times$ & $\times$ & $\times$ & \tabincell{c}{Mixed-integer linear\\programming (MILP)} \\
			\tabincell{c}{Gomez-de-Arteche-\\Botas 2025 \cite{gomez2025heat}} & Multiple 1-in-1 & $\times$ & $\times$ & $\times$ & $\times$ & $\checkmark$ & Heat pump & Mathematical method \\
			Amin 2026 \cite{amin2026modelling} & 1-in-1 & $\times$ & $\checkmark$ & $\times$ & $\checkmark$ & $\checkmark$ & $\times$ & Thermal modeling \\
			Allan 2026 \cite{allan2026optimization} & Multiple 1-in-1 & $\times$ & $\checkmark$ & $\times$ & $\times$ & $\checkmark$ & DHW & \tabincell{c}{Dynamic optimized\\control} \\
			Ding 2024 \cite{ding2024study} & Multiple 1-in-1 & $\times$ & $\times$ & $\checkmark$ & $\times$ & $\times$ & $\times$ & MILP scheduling \\
			Zhong 2026 \cite{zhong2026real} & 1-in-1 & $\times$ & $\checkmark$ & $\times$ & $\times$ & $\times$ & $\times$ & Real-time MPC \\
			Han 2024 \cite{han2024dual} & Multiple 1-in-1 & $\times$ & $\times$ & $\checkmark$ & $\checkmark$ & $\checkmark$ & DHN & \tabincell{c}{Dual-layer model\\predictive control} \\
			Han 2025 \cite{han2025robust} & Multiple 1-in-1 & $\times$ & $\times$ & $\checkmark$ & $\checkmark$ & $\checkmark$ & DHN & \tabincell{c}{Two-stage robust\\scheduling} \\
			\hline
			\textbf{This work} & \tabincell{c}{Large-scale plant with\\multiple $N$-in-1 groups} & $\checkmark$ & $\checkmark$ & $\checkmark$ & $\checkmark$ & $\checkmark$ & PTUS & \tabincell{c}{Hierarchical electricity-heat\\coordination framework} \\
			\hline\hline
		\end{tabular}
	\end{adjustbox}

	\vspace{6pt}
	\parbox{0.98\textwidth}{\footnotesize\emph{Note:} $\checkmark$ indicates that the process is considered, whereas $\times$ indicates that it is not considered. DHW, domestic hot water; DHN, district heating network; MPC, model predictive control.}
\end{table}

\subsection{Literature Review}
\label{sec:lit_review}

Research on thermal and energy management of hydrogen plants can be grouped into three areas: thermal modeling of AWE systems, coordinated operation of multiple stacks, and electricity-heat integration with waste heat recovery.

\emph{a) Thermal dynamics modeling of AWE systems.}
Dynamic AWE models describe the electrochemical and thermal responses of electrolyzers under fluctuating power. Existing studies have related current density and stack temperature to hydrogen production and energy efficiency \cite{jin2025alkaline, meng2026model}. Zhong et al. \cite{zhong2025improving} further analyzed variable-load operation using an unsteady thermodynamic model. These studies show that stack temperature is determined by electrochemical operation and, in turn, affects voltage efficiency, startup, and operating stability.

Most available models, however, describe a single stack or simplified AWE system. Thermal interactions among multiple stacks, BoP components, the plant building, and the PTUS are generally neglected or represented by fixed boundary conditions. They therefore cannot describe plant-wide heat generation, transfer, recovery, and dissipation.

\emph{b) Coordinated scheduling of multiple stacks.}
As ReP2H plants increase in scale, coordinated multi-stack operation has been studied to improve renewable utilization, operating economy, and flexibility. Recent work has addressed dynamic modeling and hierarchical scheduling \cite{guan2025dynamic, guan2026region}, coordinated stack operation and rolling power allocation \cite{zou2025control, wang2025collaborative}, and mean-field control for large stack clusters \cite{chen2025mean}. Startup and shutdown dynamics, thermal constraints, and safety limits have also been incorporated into scheduling models \cite{xu2025optimization, ma2025cold, qiu2023extended, firdous2025utility}.

In these studies, however, temperature is mainly treated as an operating constraint. Heat is rarely managed as a plant-wide resource that can be stored, transferred, and reused. Preheating, thermal standby, and waste heat recovery are therefore weakly coupled with stack commitment and power allocation.

\emph{c) Electricity-heat integration and waste heat recovery.}
Electricity-hydrogen-heat coordination has also been studied in integrated energy systems (IESs). Gomez-de-Arteche-Botas et al. \cite{gomez2025heat} investigated heat pumps for recovering waste heat from green hydrogen production, while Amin et al. \cite{amin2026modelling} evaluated the energy-saving potential of AWE waste heat. Han et al. developed bidirectional heat exchange and scheduling methods for AWE systems coupled to district heating networks (DHNs), allowing recovered stack heat to support external heating and external heat to assist electrolyzer thermal regulation \cite{han2024dual, han2025robust}. Related studies have considered electricity-hydrogen-heat coordination and heating-network integration \cite{allan2026optimization, ding2024study, zhong2026real, li2018operation}.

These studies demonstrate the value of bidirectional heat exchange, but mainly consider coordination between electrolyzers and external heating networks. Internal heat coordination in large hydrogen plants has received less attention. In particular, the coupled thermal dynamics of multiple $N$-in-1 groups, shared BoP, PTUS, and plant buildings, as well as fast coordination of stack preheating, thermal standby, and waste heat recovery, remain insufficiently addressed.

Table~\ref{tab:literature_review} summarizes the main differences among recent studies. Two gaps are evident. First, existing models do not fully represent plant-wide thermal coupling among the stack cluster, BoP, PTUS, and building. Second, existing energy management methods rarely coordinate minute-scale production scheduling with second-scale thermal regulation while allowing heat to flow bidirectionally between electrolysis and plant thermal systems. As a result, the effects of plant-wide thermal coordination on startup, thermal stress, degradation, and long-term energy performance remain unclear.

\subsection{Contributions of This Work}
\label{sec:contributions}

To address these gaps, this paper proposes a hierarchical electricity-heat coordination framework for large-scale ReP2H plants in cold regions. The upper layer performs rolling unit commitment and power allocation at a $15\ \text{min}$ resolution, while the lower layer coordinates the stack cluster, BoP, PTUS, and plant building on a second-scale basis. The main contributions are as follows:

\begin{enumerate}
	\item A plant-wide thermal topology is proposed to couple the AWE stack cluster, BoP, PTUS, and plant building. The corresponding dynamic model describes heat generation, transfer, recovery, and dissipation throughout the plant, including bidirectional heat exchange between the stack cluster and PTUS.
	
	\item A hierarchical electricity-heat coordination method is developed to link production scheduling with fast thermal regulation. Upper-layer MILP scheduling determines stack commitment and power allocation, while lower-layer adaptive current and thermal control coordinate preheating, thermal standby, and waste heat recovery under fluctuating renewable power.
	
	\item A plant-scale numerical study based on engineering data from an $80~\mathrm{MW}$ hydrogen project evaluates the proposed framework from second-scale startup transients to seasonal and annual operation. The effects on startup performance, hydrogen yield, energy and exergy efficiencies, degradation, and LCOH are quantified.
\end{enumerate}

The remainder of this paper is organized as follows. Section~\ref{sec:system_architecture} presents the plant-wide thermal architecture and dynamic models. Section~\ref{sec:two_layer_control} introduces the hierarchical electricity-heat coordination framework. Section~\ref{sec:results_and_discussion} presents the comparative case studies and seasonal operation results. Section~\ref{sec:conclusion} concludes the paper.

\section{Plant-Wide Thermal Architecture and Dynamic Models}
\label{sec:system_architecture}

\subsection{Integrated Plant Thermal Management Topology}
\label{subsec:plant_topology}

\begin{figure}[!tb]
	\centering
\includegraphics[scale=0.625]{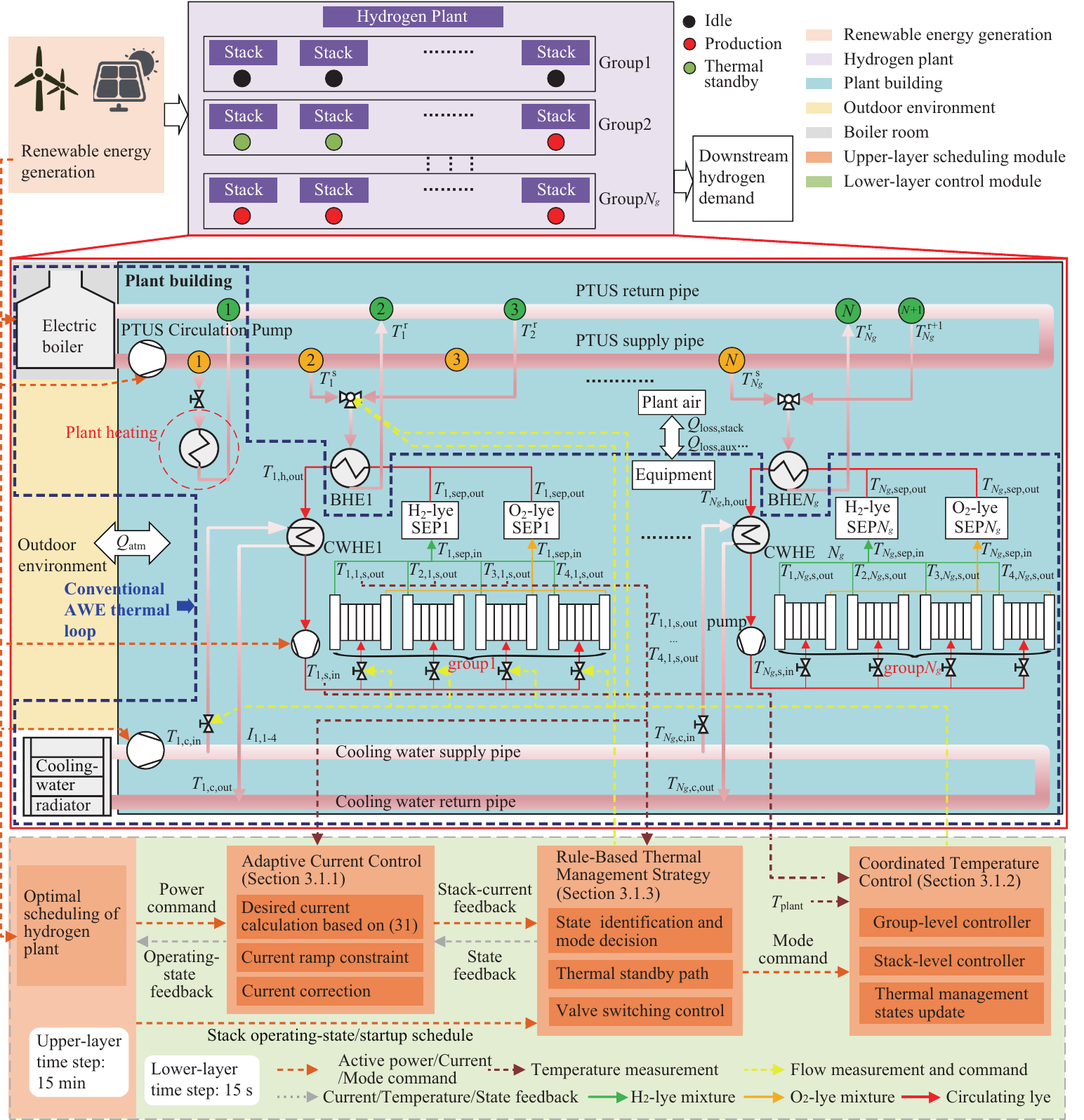}
	\vspace{-6pt}
	\caption{Proposed thermal management topology and control flowchart of the hydrogen plant.}
	\label{fig:topology_flowchart}
\end{figure}

Fig.~\ref{fig:topology_flowchart} compares the proposed hydrogen plant topology with a conventional design. The plant consists of AWE stack clusters, shared balance of plant (BoP) equipment, and a PTUS. Rectifiers supply DC power to parallel stacks. Gas-lye mixtures from the stacks enter shared gas-lye separators (SEPs), where the gas and lye are separated. The lye is then cooled by a cooling water heat exchanger (CWHE), replenished with deionized water, and recirculated to the stacks. Each group shares one CWHE connected to the plant cooling water loop \cite{qiu2026dynamic}.

In the conventional topology, enclosed by the blue dashed lines in Fig.~\ref{fig:topology_flowchart}, the AWE thermal loop and PTUS operate independently. Stack heat cannot be transferred to the PTUS, and PTUS heat cannot be used to regulate the lye temperature. To enable bidirectional thermal coupling, this work introduces a bidirectional heat exchanger (BHE) upstream of each CWHE and connects it to the PTUS.

\begin{figure}[tb]
	\centering
	\includegraphics{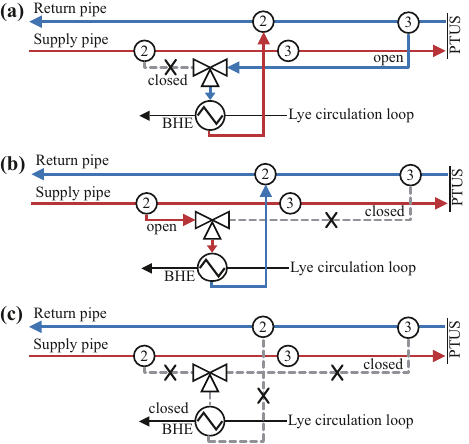}
	\caption{Three operating modes of the BHE. (a) Waste heat recovery; (b)  Preheating/thermal standby; and (c) Bypass.}
	\label{fig:three_modes}
\end{figure}

At Node~2, for example, the BHE is connected to the PTUS through a thermal reversing valve and a bypass branch linked to the downstream return node (Node~3). The valve switches the BHE among the three modes shown in Fig.~\ref{fig:three_modes} according to the thermal state and power input.
\begin{itemize}
	\item \textbf{Recovery mode:} The BHE is connected to the PTUS return pipe at Node~3 and isolated from the supply pipe at Node~2. The return water absorbs stack waste heat and transfers it to the PTUS.
	
	\item \textbf{Preheating/thermal standby mode:} The BHE is connected to the PTUS supply pipe at Node~2 and isolated from the return pipe at Node~3. High-temperature supply water heats the circulating lye for stack preheating or thermal standby.
	
	\item \textbf{Bypass mode:} Both PTUS connections are closed, isolating the lye loop from the PTUS during full shutdown or when heat exchange is unnecessary.
\end{itemize}

The proposed topology therefore allows stack waste heat to support plant heating while enabling the PTUS to provide heat for startup and thermal standby. Heat released by the stacks, BoP, and PTUS also affects the plant air temperature through equipment dissipation, radiators, envelope heat transfer, and ventilation. These heat paths couple the stack cluster, BoP, PTUS, and plant building within one thermal system.

\subsection{Mass and Heat Transfer Model of the AWE System}

As shown in Fig.~\ref{fig:topology_flowchart}, multiple $N$-in-1 AWE groups are thermally coupled to the PTUS through shared BoP components, including SEPs, BHEs, and CWHEs. The electrochemical and thermal models are based on \cite{qi2023thermal, qiu2026dynamic}, with additional equations introduced for the shared BoP and PTUS coupling. The governing equations are summarized in Table~\ref{tab:model_summary}. The subscript $(i,j)$ denotes stack $i$ in group $j$, and $T_{j,\mathrm{s,in}}$ is the common inlet temperature of all stacks in group $j$.

Specifically, \eqref{eq:UItotal1} describes the stack voltage, while \eqref{eq:faraday2} accounts for stray-current effects in the lye channels when calculating Faradaic efficiency. Hydrogen production and electrolytic power are given by \eqref{eq:hydrogenflow3} and \eqref{eq:power4}, respectively. The thermal dynamics of the stack, SEP, BHE, and CWHE are described by \eqref{eq:stack5}--\eqref{eq:stack7}, \eqref{eq:sep9}--\eqref{eq:sep11}, \eqref{eq:he12}--\eqref{eq:he14}, and \eqref{eq:cwhe15}--\eqref{eq:cwhe17}, respectively. For the BHE, $b_{\mathrm{hc}}=1$ denotes PTUS-supplied preheating or thermal standby, whereas $b_{\mathrm{hc}}=0$ denotes waste heat recovery to the PTUS return loop. In bypass mode, $v_{j,\mathrm{pre}}=v_{j,\mathrm{rec}}=0$, and the BHE heat transfer rate is zero. Heat-transfer and electrochemical parameters can also be updated through online identification \cite{qiu2023dynamic_parameter}.

\label{subsec:mass_heat_model}

\begin{table}[!tbp]\scriptsize
	\renewcommand{\arraystretch}{1.0}
	\setlength{\tabcolsep}{5pt}
	\setlength{\abovedisplayskip}{-6pt}
	\setlength{\belowdisplayskip}{-3pt}
	\setlength{\abovedisplayshortskip}{2pt}
	\setlength{\belowdisplayshortskip}{2pt}
	\caption{Summary of the AWE thermal dynamics, hydrogen production, and plant building heat transfer models.}
	\vspace{6pt}
	\label{tab:model_summary}
	\centering
	
	\begin{tabular}{ >{\centering\arraybackslash}m{2.4cm} >{\raggedright\arraybackslash}m{13.3cm} }
		\hline\hline
		Submodel & \multicolumn{1}{c}{Physical Process Model} \\
		\hline
		
		\tabincell{c}{Electrochemistry\\and Production}
		&
		{\begin{flalign}
				& U_{i,j}^{\mathrm{cell}}
				= U^{\mathrm{rev}}
				+ (r_1+r_2T_{i,j,\mathrm{s,out}}+r_3p)I_{i,j}
				+ s\log\left[\left(t_1+\frac{t_2}{T_{i,j,\mathrm{s,out}}}
				+\frac{t_3}{T_{i,j,\mathrm{s,out}}^2}\right)I_{i,j}+1\right] && \label{eq:UItotal1}\\
				& \eta_{i,j}^{\mathrm{cell}}
				{= \frac{(0.1I_{i,j})^{2}}{f_1+(0.1I_{i,j})^{2}}f_2}
				&& \label{eq:faraday2}\\
				& \dot{n}_{i,j,\mathrm{H_2}}
				= \eta_{i,j}^{\mathrm{cell}}N^{\mathrm{cell}}I_{i,j}/(2F)
				&& \label{eq:hydrogenflow3}\\
				& P_{i,j}^{\mathrm{ele}}
				= N^{\mathrm{cell}}U_{i,j}^{\mathrm{cell}}I_{i,j}
				&& \label{eq:power4}
		\end{flalign}}
		\\[-2ex]\hline
		
		\tabincell{c}{Stack Heat\\Transfer}
		&
		{\begin{flalign}
				& C_{i,j,\mathrm{s}}
				\frac{\mathrm{d}T_{i,j,\mathrm{s,out}}}{\mathrm{d}t}
				= Q_{i,j,\mathrm{ele}}
				- Q_{i,j,\mathrm{s,diss}}
				- c_{\mathrm{lye}}v_{i,j,\mathrm{lye}}\rho_{\mathrm{lye}}
				(T_{i,j,\mathrm{s,out}}-T_{j,\mathrm{s,in}})
				&& \label{eq:stack5}\\
				& Q_{i,j,\mathrm{ele}}
				= \eta_{i,j}^{\mathrm{cell}}N^{\mathrm{cell}}I_{i,j}
				(U_{i,j}^{\mathrm{cell}}-U^{\mathrm{th}})
				+ (1-\eta_{i,j}^{\mathrm{cell}})N^{\mathrm{cell}}I_{i,j}U_{i,j}^{\mathrm{cell}}
				&& \label{eq:stack6}\\
				& Q_{i,j,\mathrm{s,diss}}
				= h_{\mathrm{s}}A_{\mathrm{s,diss}}(T_{i,j,\mathrm{s,out}}-T_{\mathrm{plant}})
				+ \sigma A_{\mathrm{s,diss}}\varepsilon_{\mathrm{s}}
				(T_{i,j,\mathrm{s,out}}^{4}-T_{\mathrm{plant}}^{4})
				&& \label{eq:stack7}\\
				& h_{\mathrm{s}}
				= 2.51 \times 0.52
				\left((T_{i,j,\mathrm{s,out}}-T_{\mathrm{plant}})/\varphi_{\mathrm{s}}\right)^{0.25}
				&& \label{eq:stack8}
		\end{flalign}}
		\\[-2ex]\hline
		
		\tabincell{c}{SEP Heat\\Transfer}
		&
		{\begin{flalign}
				& C_{j,\mathrm{sep}}
				\frac{\mathrm{d}T_{j,\mathrm{sep,out}}}{\mathrm{d}t}
				= \frac{1}{2}c_{\mathrm{lye}}v_{j,\mathrm{tot,lye}}\rho_{\mathrm{lye}}
				(T_{j,\mathrm{sep,in}}-T_{j,\mathrm{sep,out}})
				- Q_{j,\mathrm{sep,diss}}
				&& \label{eq:sep9}\\
				& {
					Q_{j,\mathrm{sep,diss}}
					= h_{\mathrm{sep}}A_{\mathrm{sep,diss}}(T_{j,\mathrm{sep,out}}-T_{\mathrm{plant}})
					+ \sigma A_{\mathrm{sep,diss}}\varepsilon_{\mathrm{sep}}
					(T_{j,\mathrm{sep,out}}^{4}-T_{\mathrm{plant}}^{4})
				}
				&& \label{eq:sep10}\\
				& T_{j,\mathrm{sep,in}}
				= \left(\sum\nolimits_{i=1}^{N}v_{i,j,\mathrm{lye}}T_{i,j,\mathrm{s,out}}\right)/v_{j,\mathrm{tot,lye}}
				&& \label{eq:sep11}
		\end{flalign}}
		\\ \hline
		
		\tabincell{c}{BHE Heat\\Transfer}
		&
		{\begin{flalign}
				& C_{j,\mathrm{he}}
				\frac{\mathrm{d}T_{j,\mathrm{h,out}}}{\mathrm{d}t}
				= c_{\mathrm{lye}}v_{j,\mathrm{tot,lye}}\rho_{\mathrm{lye}}
				(T_{j,\mathrm{sep,out}}-T_{j,\mathrm{h,out}})
				- (T_{j,\mathrm{h,out}}-T_{\mathrm{plant}})/{R_{\mathrm{BHE}}}
				{- k_{\mathrm{h}}A_{\mathrm{h}}\Delta T_{j,\mathrm{h}}}
				&& \label{eq:he12}\\
				& C_{j,\mathrm{h}}
				\frac{\mathrm{d}T_{j}^{\mathrm{HD,r}}}{\mathrm{d}t}
				= c_{\mathrm{h}}\rho_{\mathrm{h}}(v_{j,\mathrm{hs}}+v_{j+1,\mathrm{hr}})
				(b_{\mathrm{hc}}T_{j}^{\mathrm{s}}+(1-b_{\mathrm{hc}})T_{j+1}^{\mathrm{r}}-T_{j}^{\mathrm{HD,r}})
				{+ k_{\mathrm{h}}A_{\mathrm{h}}\Delta T_{j,\mathrm{h}}}
				&& \label{eq:he13}\\
				& \Delta T_{j,\mathrm{h}}
				= \frac{(T_{j,\mathrm{h,out}}-T_{j}^{\mathrm{HD,r}})
					-(T_{j,\mathrm{sep,out}}-b_{\mathrm{hc}}T_{j}^{\mathrm{s}}-(1-b_{\mathrm{hc}})T_{j+1}^{\mathrm{r}})}
				{\log\left[(T_{j,\mathrm{h,out}}-T_{j}^{\mathrm{HD,r}})/
					(T_{j,\mathrm{sep,out}}-b_{\mathrm{hc}}T_{j}^{\mathrm{s}}-(1-b_{\mathrm{hc}})T_{j+1}^{\mathrm{r}})\right]}
				&& \label{eq:he14}
		\end{flalign}}
		\\[2ex]\hline
		
		\tabincell{c}{CWHE Heat\\Transfer}
		&
		{\begin{flalign}
				& C_{j,\mathrm{ce}}
				\frac{\mathrm{d}T_{j,\mathrm{s,in}}}{\mathrm{d}t}
				= c_{\mathrm{lye}}v_{j,\mathrm{tot,lye}}\rho_{\mathrm{lye}}
				(T_{j,\mathrm{h,out}}-T_{j,\mathrm{s,in}})
				- (T_{j,\mathrm{s,in}}-T_{\mathrm{plant}})/{R_{\mathrm{CWHE}}}
				{- k_{\mathrm{c}}A_{\mathrm{c}}\Delta T_{j,\mathrm{c}}}
				&& \label{eq:cwhe15}\\
				& C_{j,\mathrm{c}}
				\frac{\mathrm{d}T_{j,\mathrm{c,out}}}{\mathrm{d}t}
				= c_{\mathrm{c}}v_{j,\mathrm{c}}\rho_{\mathrm{c}}
				(T_{j,\mathrm{c,in}}-T_{j,\mathrm{c,out}})
				{+ k_{\mathrm{c}}A_{\mathrm{c}}\Delta T_{j,\mathrm{c}}}
				&& \label{eq:cwhe16}\\
				& \Delta T_{j,\mathrm{c}}
				= \frac{(T_{j,\mathrm{s,in}}-T_{j,\mathrm{c,out}})
					-(T_{j,\mathrm{h,out}}-T_{j,\mathrm{c,in}})}
				{\log\left[(T_{j,\mathrm{s,in}}-T_{j,\mathrm{c,out}})/
					(T_{j,\mathrm{h,out}}-T_{j,\mathrm{c,in}})\right]}
				&& \label{eq:cwhe17}
		\end{flalign}}
		\\[0ex]\hline
		
		\tabincell{c}{Plant Building\\Thermal Balance}
		&
		{\begin{flalign}
				& C_{\mathrm{plant}}\frac{\mathrm{d}T_{\mathrm{plant}}}{\mathrm{d}t}
				= N_{\mathrm{rad}}U_{\mathrm{rad}}A_{\mathrm{rad}}
				(T_{\mathrm{rad}}-T_{\mathrm{plant}})
				+ \sum\nolimits_{j=1}^{N_{\mathrm{g}}}\sum\nolimits_{i=1}^{N}Q_{i,j,\mathrm{s,diss}} && \nonumber \\
				& \quad + \sum\nolimits_{j=1}^{N_{\mathrm{g}}}
				\left(
				Q_{j,\mathrm{sep,diss}}
				+ (T_{j,\mathrm{h,out}}-T_{\mathrm{plant}})/{R_{\mathrm{BHE}}}
				+ (T_{j,\mathrm{s,in}}-T_{\mathrm{plant}})/{R_{\mathrm{CWHE}}}
				\right)
				- Q_{\mathrm{atm}}-Q_{\mathrm{vent}} \label{eq:plant_thermal_balance}\\
				& Q_{\mathrm{atm}}
				= U_{\mathrm{b}}A_{\mathrm{b}}(T_{\mathrm{plant}}-T_{\mathrm{atm}}),
				\quad
				Q_{\mathrm{vent}}
				= \alpha_{\mathrm{vent}}U_{\mathrm{vent,max}}
				(T_{\mathrm{plant}}-T_{\mathrm{atm}})
				&& \label{eq:plant_heat_loss}
		\end{flalign}}
		\\[-2ex]\hline
		
		\tabincell{c}{Radiator Heat\\Transfer}
		&
		{\begin{flalign}
				& C_{\mathrm{rad}}
				\frac{\mathrm{d}T_{\mathrm{rad}}}{\mathrm{d}t}
				= N_{\mathrm{rad}}U_{\mathrm{water}}A_{\mathrm{water}}
				\left((T_{\mathrm{w}}^{\mathrm{s}}+T_{\mathrm{w}}^{\mathrm{r}})/2
				-T_{\mathrm{rad}}\right)
				- N_{\mathrm{rad}}U_{\mathrm{rad}}A_{\mathrm{rad}}
				(T_{\mathrm{rad}}-T_{\mathrm{plant}})
				&& \label{eq:radiator_temperature_dynamic}\\
				& U_{\mathrm{rad}}
				= {1.31\left|T_{\mathrm{rad}}-T_{\mathrm{plant}}\right|^{1/3}
					+ \varepsilon_{\mathrm{rad}}\sigma (T_{\mathrm{rad}}+T_{\mathrm{plant}})
					(T_{\mathrm{rad}}^{2}+T_{\mathrm{plant}}^{2})}
				&& \label{eq:radiator_heat_transfer_coefficient}\\
				& C_{\mathrm{w}}
				\frac{\mathrm{d}T_{\mathrm{w}}^{\mathrm{r}}}{\mathrm{d}t}
				= c_{\mathrm{h}}v_{\mathrm{h}}\rho_{\mathrm{h}}
				(T_{\mathrm{w}}^{\mathrm{s}}-T_{\mathrm{w}}^{\mathrm{r}})
				-N_{\mathrm{rad}} U_{\mathrm{water}}A_{\mathrm{water}}
				\left((T_{\mathrm{w}}^{\mathrm{s}}+T_{\mathrm{w}}^{\mathrm{r}})/2
				-T_{\mathrm{rad}}\right)
				&& \label{eq:radiator_thermal_dynamic}
		\end{flalign}}
		\\[-1ex]
		
		\hline\hline
	\end{tabular}
\end{table}

\subsection{Thermal Dynamics of the PTUS and the Plant Building}
\label{subsec:factory_building_model}

\begin{table}[tb]\scriptsize
	\renewcommand{\arraystretch}{1.0}
	\setlength{\tabcolsep}{5pt}
	\setlength{\abovedisplayskip}{-6pt}
	\setlength{\belowdisplayskip}{-2pt}
	\setlength{\abovedisplayshortskip}{2pt}
	\setlength{\belowdisplayshortskip}{2pt}
	\caption{Summary of the hydraulic and thermal models of the PTUS.}
	\vspace{6pt}
	\label{tab:hps_model_summary}
	\centering
	
	\begin{tabular}{>{\centering\arraybackslash}m{2.4cm} >{\raggedright\arraybackslash}m{13.3cm}}
		\hline\hline
		Submodel & \multicolumn{1}{c}{Physical Process Model} \\
		\hline
		
		\tabincell{c}{Flow\\Balance}
		&
		{\begin{flalign}
				& v_{j}^{\mathrm{HS}}+\sum\nolimits_{k\in\mathcal{I}(j)}v_{kj}
				= v_{j}^{\mathrm{HD}}+\sum\nolimits_{k\in\mathcal{O}(j)}v_{jk}
				&& \label{eq:hydraulic_balance}\\
				& p_{\mathrm{w},j}-p_{\mathrm{w},k}=6.88\times10^{-3}{K^{0.25}l}/({\rho_{\mathrm{h}}d^{1.25}})\times v_{jk}^{2} && \label{eq:hps_pressure_drop}\\
				& P_{\mathrm{pump}}={v_{\mathrm{pump}}\rho_{\mathrm{h}}gH_{\mathrm{p}}}/{\eta_{\mathrm{pump}}}
				&& \label{eq:HP_power27}
		\end{flalign}}
		\\[-2ex]
		\hline
		
		\tabincell{c}{Electric Boiler\\Heat Supply}
		&
		{\begin{flalign}
				& P_{\mathrm{boiler}}^{\mathrm{out}}=\eta_{\mathrm{boiler}}P_{\mathrm{boiler}}^{\mathrm{in}} && \label{eq:boiler_efficiency}\\
				& P_{\mathrm{boiler}}^{\mathrm{out}}=c_{\mathrm{h}}\rho_{\mathrm{h}}v_{j}^{\mathrm{HS}}\left(T_{j}^{\mathrm{HS,s}}-T_{j}^{\mathrm{HS,r}}\right) && \label{eq:boiler_power29}
		\end{flalign}}
		\\[-2ex]
		\hline
		
		\tabincell{c}{Pipe Heat\\Transfer}
		&
		{\begin{flalign}
				& Q_{jk}=c_{\mathrm{h}}\rho_{\mathrm{h}}v_{jk}
				\left(1-\psi_{jk}^{\mathrm{HL}}\right)
				\left(T_{j}^{\mathrm{s}}-T_{\mathrm{plant}}\right)
				&& \label{eq:hps_pipe_heat_transfer}
		\end{flalign}}
		\\[-2ex]
		\hline
		
		\tabincell{c}{Nodal Temperature\\Balance}
		&
		{\begin{flalign}
				& (T_{j}^{\mathrm{s}}-T_{\mathrm{plant}})\left(v_{j}^{\mathrm{HS}}+\sum\nolimits_{k\in\mathcal{I}(j)}v_{kj}\right)=(T_{j}^{\mathrm{HS,s}}-T_{\mathrm{plant}})v_{j}^{\mathrm{HS}}+\sum\nolimits_{k\in\mathcal{I}(j)}\psi_{kj}^{\mathrm{HL}}(T_{k}^{\mathrm{s}}-T_{\mathrm{plant}})v_{kj} && \label{eq:hps_supply_temperature}\\
				& (T_{j}^{\mathrm{r}}-T_{\mathrm{plant}})\left(v_{j}^{\mathrm{HD}}+\sum\nolimits_{k\in\mathcal{O}(j)}v_{jk}\right)=(T_{j}^{\mathrm{HD,r}}-T_{\mathrm{plant}})v_{j}^{\mathrm{HD}}+\sum\nolimits_{k\in\mathcal{O}(j)}\psi_{jk}^{\mathrm{HL}}(T_{k}^{\mathrm{r}}-T_{\mathrm{plant}})v_{jk} && \label{eq:hps_return_temperature}
		\end{flalign}}
		\\[-2ex]
		
		\hline\hline
	\end{tabular}
\end{table}

Table~\ref{tab:model_summary} also summarizes the thermal models of the plant building and radiators. Heat dissipated by the stack cluster, shared BoP, and PTUS enters the plant air, while radiators provide additional heating. Heat is lost to the outdoor environment through the building envelope and ventilation.

The plant air temperature $T_{\mathrm{plant}}$ is represented by a spatially uniform lumped thermal node whose equivalent thermal capacity includes the indoor air and the effective thermal mass of the building and installed equipment \cite{chi2023hvac,tol2023development}. Its dynamics are given by \eqref{eq:plant_thermal_balance} and \eqref{eq:plant_heat_loss}. Equations~\eqref{eq:radiator_temperature_dynamic}--\eqref{eq:radiator_thermal_dynamic} describe heat transfer from the circulating water to the plant air by convection and radiation \cite{liu2024radiator}.

The hydraulic and thermal models of the PTUS are summarized in Table~\ref{tab:hps_model_summary}. In \eqref{eq:hydraulic_balance}, $\mathcal{I}(j)$ and $\mathcal{O}(j)$ denote the upstream and downstream nodes connected to node $j$, respectively. The branch temperature retention factors $\psi_{jk}^{\mathrm{HL}}$ are specified in Table~\ref{tab:ptus_parameters}. The flow balance, pipe heat loss, and nodal temperature mixing equations follow established hydraulic and thermal network models \cite{shabanpour2015integrated,ahmed2014strategic}. The common hydraulic and auxiliary parameters used in the case study, including $K$, $H_{\mathrm{p}}$, $\eta_{\mathrm{pump}}$, $\eta_{\mathrm{boiler}}$, and $T_{1}^{\mathrm{HS,s}}$, are listed in Table~\ref{tab:plant_parameters}.

\section{Hierarchical Electricity-Heat Coordination Framework for Hydrogen Plants}
\label{sec:two_layer_control}

\begin{figure}[!htb]
	\centering
\includegraphics[scale=0.75]{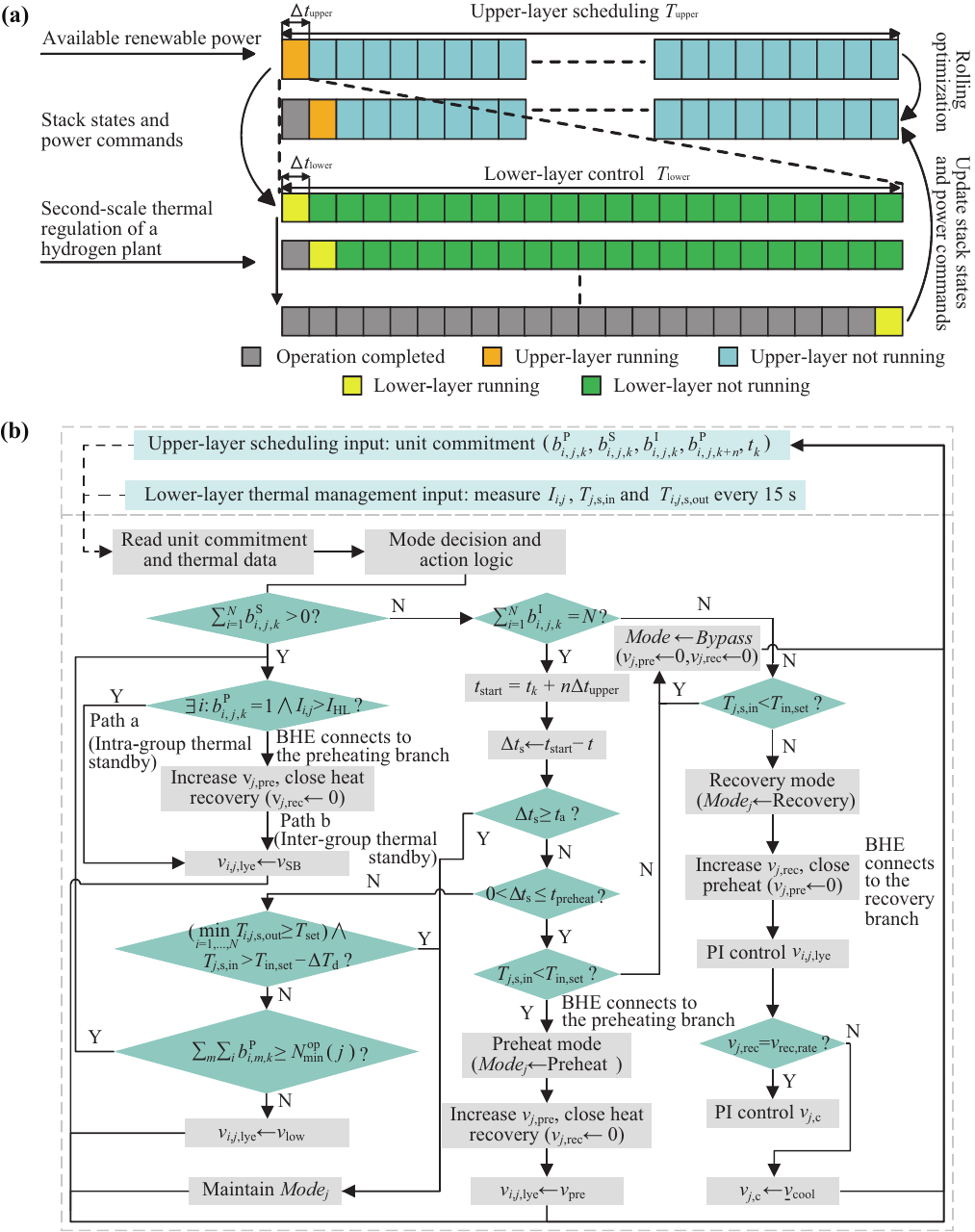}
	\vspace{-6pt}
	\caption{(a) Schematic of the hierarchical electricity-heat coordination framework. (b) Thermal management strategy for stack operating states and heat exchange modes.}
	\label{fig:bilevel_control}
\end{figure}

The hierarchical electricity-heat coordination framework uses the integrated thermal model in Section~\ref{sec:system_architecture}. As shown in Fig.~\ref{fig:bilevel_control}(a), the upper layer performs rolling unit commitment and stack power allocation over horizon $T_{\mathrm{upper}}$ with step $\Delta t_{\mathrm{upper}}$ based on renewable power forecasts. The lower layer operates over period $T_{\mathrm{lower}}$ with step $\Delta t_{\mathrm{lower}}=15~\mathrm{s}$ and executes these commands through adaptive current control, coordinated temperature control, and the thermal management strategy described in Section~\ref{subsec:lower_layer_control}.

At each lower-layer step, stack currents and heat-exchange modes are adjusted subject to voltage, temperature, and plant thermal constraints. At the end of each upper-layer interval, the executed stack power and measured thermal states are fed back to update the next rolling schedule.

\subsection{Lower-Layer Thermal Management of the Hydrogen Plant at Second-Scale Resolution}
\label{subsec:lower_layer_control}
\subsubsection{Adaptive Current Control}
\label{subsubsec:adaptive_current_control}

{Within upper-layer scheduling interval $k$, the lower layer updates the stack current and thermal states at control step $q$ with step size $\Delta t_{\mathrm{lower}}$. Given the scheduled stack electrochemical power command $P_{i,j,k}^{\mathrm{ele}}$, the desired current is calculated directly from this power command and the cell voltage at the previous step:}
\begin{equation}
  {
  I_{i,j,q}^{\mathrm{des}}
  =
  {P_{i,j,k}^{\mathrm{ele}}}/(
  {N^{\mathrm{cell}}U_{i,j,q-1}^{\mathrm{cell}}}),
  }
  \label{eq:desired_current}
\end{equation}

The maximum allowable current follows from the inverse voltage-current relationship:
\begin{equation}
 I_{i,j,q}^{\mathrm{max}} = f_U^{-1}(U^{\mathrm{max}}, T_{i,j,q,\mathrm{s,out}}),
\end{equation}
where $f_U^{-1}(\cdot)$ is the inverse voltage model at the current stack temperature. The current is subject to the ramping constraint
\begin{equation}
  |I_{i,j,q} - I_{i,j,q-1}| \le \Delta I_{i,j}^{\mathrm{max}} \Delta t_{\mathrm{lower}},
\end{equation}
where $\Delta I_{i,j}^{\mathrm{max}}$ is the maximum current ramp rate.

Combining the power command, voltage limit, and ramping constraint gives the executable current:
\begin{equation}
    I_{i,j,q}
    =
    \min\left\{
    I_{i,j,q}^{\mathrm{max}},
    \max\left[
    I_{i,j,q-1}-\Delta I_{i,j}^{\mathrm{max}}\Delta t_{\mathrm{lower}},
    \min\left(
    I_{i,j,q-1}+\Delta I_{i,j}^{\mathrm{max}}\Delta t_{\mathrm{lower}},
    I_{i,j,q}^{\mathrm{des}}
    \right)
    \right]
    \right\}.
\end{equation}


\subsubsection{Coordinated Temperature Control}

Because stacks within an $N$-in-1 group share the BoP and lye circulation loop, their temperatures are coupled. A two-level PI structure is therefore used to regulate both group inlet and individual stack temperatures. The group-level gains $K_{\mathrm{p,g}}$ and $K_{\mathrm{i,g}}$ and stack-level gains $K_{\mathrm{p,s}}$ and $K_{\mathrm{i,s}}$ are initialized from step-response tests and refined through simulation to balance response speed and stability while limiting integral saturation \cite{qi2023design}.

\begin{itemize}
	\item \textbf{Group level:} The PTUS water flow is regulated to maintain the common group inlet temperature $T_{j,\mathrm{s,in}}$ at $T_{\mathrm{in,set}}$ via a PI controller.
	
	\item \textbf{Stack level:} Each stack's lye flow is regulated independently to maintain $T_{i,j,\mathrm{s,out}}$ near its target temperature with an independent PI controller.
\end{itemize}


\subsubsection{Thermal Management Strategy for Stack Operating States and Heat-Exchange Modes}
\label{subsubsec:rule_based_thermal_control}

The strategy in Fig.~\ref{fig:bilevel_control}(b) maps the upper-layer production, thermal standby, and idle commands, together with lower-layer temperature measurements, to BHE operating modes and flow commands. The control logic is applicable to general $N$-in-1 configurations, whereas its thresholds, flow settings, heat-availability criterion, and PI gains depend on the plant. For the $80~\mathrm{MW}$ plant in the case study, these parameters are determined from field limits and offline thermal simulations and are listed in Table~\ref{tab:thermal_management_parameters} in \ref{sec:appendix_parameters}.

\paragraph{Rule 1: Stack state identification and thermal management mode decision}

For group $j$, the controller determines $\mathrm{Mode}j$ from the stack currents $I{i,j}$, group inlet temperature $T_{j,\mathrm{s,in}}$, stack temperatures $T_{i,j,\mathrm{s,out}}$, and upper-layer commands $b_{i,j,k}^{\mathrm{P}}$, $b_{i,j,k}^{\mathrm{S}}$, and $b_{i,j,k}^{\mathrm{I}}$. If
$\sum_{i=1}^{N}b_{i,j,k}^{\mathrm{S}}>0$, Rule~2 selects the thermal standby path. Otherwise, the controller determines whether the entire group is idle.

If
$\sum_{i=1}^{N}b_{i,j,k}^{\mathrm{I}}=N$, let $n$ be the number of scheduling intervals before the next production command satisfying
$\sum_{i=1}^{N}b_{i,j,k+n}^{\mathrm{P}}>0$. The scheduled startup time and remaining time before startup are
\[
t_{\mathrm{start}}=t_k+n\Delta t_{\mathrm{upper}},
\qquad
\Delta t_{\mathrm{s}}=t_{\mathrm{start}}-t.
\]

The idle-state logic is then defined as follows:
\begin{enumerate}
	
	\item If $\Delta t_{\mathrm{s}}\geq t_{\mathrm{a}}$, the current thermal management mode is maintained.
	
	\item
	
	If $t_{\mathrm{preheat}}<\Delta t_{\mathrm{s}}<t_{\mathrm{a}}$, the current mode is maintained when $\min_{i=1,\ldots,N}T_{i,j,\mathrm{s,out}}\geq T_{\mathrm{set}}$ and $T_{j,\mathrm{s,in}}>T_{\mathrm{in,set}}-\Delta T_{\mathrm{d}}$. Otherwise, plant-wide waste heat availability is assessed. The minimum number of production stacks required to support thermal standby is determined by
	\begin{equation}
		N_{\min}^{\mathrm{op}}(j)
		=
		\min
		\left\{
		{n_{\mathrm{op}}}\in\mathbb{Z}_{+}:
		\underline{\dot{Q}}_{\mathrm{av}}({n_{\mathrm{op}}})
		\geq
		\dot{Q}_{j,\mathrm{SB}}
		+
		\dot{Q}_{j,\mathrm{loss}}
		\right\},
		\label{eq:offline_heat_threshold}
	\end{equation}
	where $\underline{\dot{Q}}_{\mathrm{av}}({n_{\mathrm{op}}})$ is the conservative lower bound of recoverable waste heat from ${n_{\mathrm{op}}}$ production stacks over the admissible load and temperature ranges, $\dot{Q}_{j,\mathrm{SB}}$ is the standby heat demand of group $j$, and $\dot{Q}_{j,\mathrm{loss}}$ represents heat-transfer and distribution losses. Waste heat is considered sufficient when
	\begin{equation}
		N_{\mathrm{op},k}
		=
		\sum\nolimits_{m=1}^{N_{\mathrm{g}}}
		\sum\nolimits_{i=1}^{N}
		b_{i,m,k}^{\mathrm{P}}
		\geq
		N_{\min}^{\mathrm{op}}(j).
		\label{eq:available_waste_heat}
	\end{equation}
	
	The threshold $N_{\min}^{\mathrm{op}}(j)$ is evaluated offline using the thermal model in Section~\ref{sec:system_architecture}. For the four-group, seven-node network in the case study, $\left[N_{\min}^{\mathrm{op}}(1),N_{\min}^{\mathrm{op}}(2),N_{\min}^{\mathrm{op}}(3),N_{\min}^{\mathrm{op}}(4)\right]=[1,\,3,\,6,\,9]$. These values depend on plant topology and operating parameters and must be recalibrated for other configurations. If \eqref{eq:available_waste_heat} is satisfied, Rule~2 is activated. Otherwise, $v_{i,j,\mathrm{lye}}\leftarrow v_{\mathrm{low}}$, and the current mode is retained.

	\item If $0<\Delta t_{\mathrm{s}}\leq t_{\mathrm{preheat}}$, Preheat mode is activated when $T_{j,\mathrm{s,in}}<T_{\mathrm{in,set}}$. The BHE is connected to the PTUS supply branch, $v_{j,\mathrm{pre}}$ is increased, $v_{j,\mathrm{rec}}\leftarrow0$, and $v_{i,j,\mathrm{lye}}\leftarrow v_{\mathrm{pre}}$. Otherwise, Bypass mode is activated with $v_{j,\mathrm{pre}}\leftarrow0$ and $v_{j,\mathrm{rec}}\leftarrow0$.
	
\end{enumerate}

If
$\sum_{i=1}^{N}b_{i,j,k}^{\mathrm{I}}\neq N$
and no stack is in thermal standby, at least one stack in group $j$ is producing hydrogen. Bypass mode is used when
$T_{j,\mathrm{s,in}}<T_{\mathrm{in,set}}$,
with
$v_{j,\mathrm{pre}}\leftarrow0$
and
$v_{j,\mathrm{rec}}\leftarrow0$.
Otherwise, Recovery mode is activated: the BHE is connected to the PTUS return branch,
$v_{j,\mathrm{rec}}$ is increased,
$v_{j,\mathrm{pre}}\leftarrow0$, and
$v_{i,j,\mathrm{lye}}$ is regulated by the PI controller.
Active cooling-water regulation is enabled only when
$v_{j,\mathrm{rec}}=v_{\mathrm{rec,rate}}$;
otherwise,
$v_{j,\mathrm{c}}\leftarrow v_{\mathrm{cool}}$.

\paragraph{Rule 2: Thermal standby path selection}
Thermal standby is supplied either within the group or through the PTUS. A production stack satisfying
$b_{i,j,k}^{\mathrm{P}}=1$
and
$I_{i,j}>I_{\mathrm{HL}}$
is classified as a high-load stack.

\begin{enumerate}
	
	
	\item \textit{Intra-group thermal standby:}
	If $\exists\,i$ such that $b_{i,j,k}^{\mathrm{P}}=1\land I_{i,j}>I_{\mathrm{HL}}$, waste heat from the high-load production stack is redistributed through the shared lye circulation loop. The lye flow of each standby stack is set to $v_{i,j,\mathrm{lye}}\leftarrow v_{\mathrm{SB}}$.
	
	\item \textit{Inter-group thermal standby:}
	If no stack satisfies $b_{i,j,k}^{\mathrm{P}}=1\land I_{i,j}>I_{\mathrm{HL}}$, local waste heat is considered insufficient. The BHE is then connected to the PTUS supply branch, $v_{j,\mathrm{pre}}$ is increased, $v_{j,\mathrm{rec}}\leftarrow0$, and the lye flow of each standby stack is set to $v_{i,j,\mathrm{lye}}\leftarrow v_{\mathrm{SB}}$.
	
\end{enumerate}

After either standby path is activated, it is maintained while
$\min_{i=1,\ldots,N}T_{i,j,\mathrm{s,out}}\geq T_{\mathrm{set}}$
and
$T_{j,\mathrm{s,in}}>T_{\mathrm{in,set}}-\Delta T_{\mathrm{d}}$.
Otherwise, \eqref{eq:available_waste_heat} is reevaluated. If sufficient waste heat is available, the controller again selects the intra-group or inter-group path according to the high-load criterion. Otherwise,
$v_{i,j,\mathrm{lye}}\leftarrow v_{\mathrm{low}}$,
and the current mode is maintained.

\subsection{Upper-Layer Unit Commitment and Power Allocation}
\label{subsec:upper_layer_scheduling}

The upper layer operates at a $15~\mathrm{min}$ resolution and determines stack commitment and power allocation. Explicitly including stack temperature dynamics can better represent thermal inertia \cite{qiu2023extended, guan2026region}, but introduces nonlinear intertemporal coupling. Because fast thermal dynamics and heat-exchange mode switching are handled by the lower layer, the upper layer retains only operating-state decisions, power allocation, and first-step thermal feasibility.

\subsubsection{Production, Thermal Standby, and Idle State Switching}
\label{subsubsec:state_switching_logic}

Each stack switches among production, thermal standby, and idle states according to renewable power availability. Because the lower-layer strategy in Section~\ref{subsubsec:rule_based_thermal_control} performs preheating before scheduled startup, the upper layer does not model the heating process explicitly. The state-switching constraints are
\begin{gather}
  b_{i,j,k}^{\mathrm{P}} + b_{i,j,k}^{\mathrm{S}} + b_{i,j,k}^{\mathrm{I}} = 1,
  \label{eq:state_exclusive} \\
  b_{i,j,k}^{\mathrm{P}} + b_{i,j,k}^{\mathrm{S}} + b_{i,j,k-1}^{\mathrm{I}} - 1 \le b_{i,j,k}^{\mathrm{SU}},
  \label{eq:startup_action} \\
  b_{i,j,k}^{\mathrm{I}} + b_{i,j,k-1}^{\mathrm{P}} + b_{i,j,k-1}^{\mathrm{S}} - 1 \le b_{i,j,k}^{\mathrm{SD}},
  \label{eq:shutdown_action} \\
  b_{i,j,k}^{\mathrm{SU}} + b_{i,j,k}^{\mathrm{SD}} \le 1.
  \label{eq:start_shutdown_exclusive}
\end{gather}

\subsubsection{Linearized Hydrogen Production Model}
\label{subsubsec:linearized_hydrogen_model}

Hydrogen production depends nonlinearly on stack power, temperature, and pressure \cite{sanchez2018semi, qi2023thermal}. Since the lower layer maintains the stack thermal state near the rated operating point over each $15~\mathrm{min}$ scheduling interval, the upper-layer model linearizes \eqref{eq:UItotal1}--\eqref{eq:power4} around this point as
\begin{equation}
  {
  F_{i,j,k}
  =
  A_{1}P_{i,j,k}^{\mathrm{ele}}
  +
  A_{2}b_{i,j,k}^{\mathrm{P}} - A_3 b_{i,j,k}^{\mathrm{SP}},
  }
  \label{eq:linear_h2_model}
\end{equation}
where $P_{i,j,k}^{\mathrm{ele}}$ is the scheduled stack electrochemical power; $A_1$ and $A_2$ characterize the linearized hydrogen-production relation, while $A_3$ represents the production loss during the transition from standby to production. The transition variable satisfies \cite{qiu2023extended, guan2026region}
\begin{equation}
  \begin{cases}
    b_{i,j,k}^{\mathrm{SP}} \le b_{i,j,k-1}^{\mathrm{S}}, \\
    b_{i,j,k}^{\mathrm{SP}} \le b_{i,j,k}^{\mathrm{P}}, \\
    b_{i,j,k}^{\mathrm{SP}} \ge b_{i,j,k-1}^{\mathrm{S}} + b_{i,j,k}^{\mathrm{P}} - 1.
  \end{cases}
  \label{eq:standby_to_production}
\end{equation}

The linearized variable $F_{i,j,k}$ allows the upper-layer problem to remain computationally efficient. After receiving $P_{i,j,k}^{\mathrm{ele}}$, the lower layer determines the executable current using Section~\ref{subsubsec:adaptive_current_control} and evaluates the actual hydrogen production rate $\dot{n}_{i,j,\mathrm{H}_2}$ from the nonlinear electrochemical model.

The renewable power balance includes electrolysis, thermal standby, pumping, and electric heating:
\begin{equation}
{
  \sum\nolimits_{j=1}^{N_{\mathrm{g}}}\sum\nolimits_{i=1}^{N}
  \left(P_{i,j,k}^{\mathrm{ele}}+P_{\mathrm{SB}}b_{i,j,k}^{\mathrm{S}}\right)
  + P_{k,\mathrm{pump}}^{(\xi)}
  + P_{k,\mathrm{heat}}^{(\xi)}
  \le P_k^{\mathrm{RE}},
}
  \label{eq:renewable_power_balance}
\end{equation}
\begin{equation}
{
  P_{\mathrm{min}} b_{i,j,k}^{\mathrm{P}}
  \le P_{i,j,k}^{\mathrm{ele}}
  \le
  P_{\mathrm{max}} b_{i,j,k}^{\mathrm{P}},
}
  \label{eq:stack_power_range}
\end{equation}

To avoid embedding the nonlinear PTUS equations in the upper-layer problem, pump and electric-boiler powers are fitted from detailed simulations based on \eqref{eq:HP_power27}--\eqref{eq:boiler_power29}. Fig.~\ref{fig:appendix_power_fit} and Table~\ref{tab:seasonal_fitting_params} in  \ref{sec:appendix_fitting}  report the fitting results. For each season $\xi$, the auxiliary powers are approximated as linear functions of total stack power:
\begin{equation}
{
  \begin{split}
    P_{k,\mathrm{pump}}^{(\xi)} &= \alpha_{\mathrm{pump}}^{(\xi)}
    \sum\nolimits_{j=1}^{N_{\mathrm{g}}} \sum\nolimits_{i=1}^{N} P_{i,j,k}^{\mathrm{ele},(\xi)}
    + b_{\mathrm{pump}}^{(\xi)}, \\
    P_{k,\mathrm{heat}}^{(\xi)} &= \alpha_{\mathrm{heat}}^{(\xi)}
    \sum\nolimits_{j=1}^{N_{\mathrm{g}}} \sum\nolimits_{i=1}^{N} P_{i,j,k}^{\mathrm{ele},(\xi)}
    + b_{\mathrm{heat}}^{(\xi)}.
  \end{split}
}
  \label{eq:seasonal_auxiliary_fitting}
\end{equation}

At each rolling update, let $k_0$ denote the current scheduling interval. First-step feasibility is enforced using the measured stack temperature and cell voltage:
\begin{equation}
{
P_{i,j,k_0}^{\mathrm{ele}}
\leq
N^{\mathrm{cell}}
U_{i,j,k_0-1}^{\mathrm{cell}}
I_{i,j}^{\max}
\left(
T_{i,j,k_0-1}^{\mathrm{actual}}
\right).
}
\label{eq:temperature_power_limit}
\end{equation}

For later intervals in the $24~\mathrm{h}$ horizon, \eqref{eq:stack_power_range} is retained because future stack temperatures are not explicitly predicted. After interval $k_0$ is executed by the $15~\mathrm{s}$ lower-layer controller, the measured stack temperature and cell voltage update \eqref{eq:temperature_power_limit} for the next rolling optimization. This receding-horizon implementation retains the MILP structure while accounting for measured thermal feasibility.

For compactness, the upper-layer exogenous/output variables, continuous decision variables, and binary variables at scheduling step $k$ are denoted as
\begin{equation}
{
  \begin{cases}
\bm{x}_k = \left( P_k^{\mathrm{RE}}, F_{i,j,k} \right), \\
    \bm{u}_k = \left( P_{i,j,k}^{\mathrm{ele}}, P_{k,\mathrm{pump}}^{(\xi)}, P_{k,\mathrm{heat}}^{(\xi)} \right), \\
    \bm{z}_k = \left( b_{i,j,k}^{\mathrm{P}}, b_{i,j,k}^{\mathrm{S}}, b_{i,j,k}^{\mathrm{I}}, b_{i,j,k}^{\mathrm{SU}}, b_{i,j,k}^{\mathrm{SD}}, b_{i,j,k}^{\mathrm{SP}} \right).
  \end{cases}
}
  \label{eq:upper_layer_variable_vectors}
\end{equation}

\subsubsection{Objective Function and Optimization Summary}
\label{subsubsec:objective_function}

The upper layer minimizes net operating cost subject to stack operating constraints and the first-step feasibility condition in \eqref{eq:temperature_power_limit}. Considering electricity cost, hydrogen revenue, auxiliary consumption, and startup and shutdown costs, the objective is set as
\begin{equation}
	{
		\begin{aligned}
			\min J
			=
			\sum_{k=1}^{N_{\mathrm{T}}}
			\Bigg\{
			&\Bigg[
			C^{\mathrm{W}}
			\left(
			\sum_{j=1}^{N_{\mathrm{g}}}\sum_{i=1}^{N}
			\left(P_{i,j,k}^{\mathrm{ele}}+P_{\mathrm{SB}}b_{i,j,k}^{\mathrm{S}}\right)
			+P_{k,\mathrm{pump}}^{(\xi)}
			+P_{k,\mathrm{heat}}^{(\xi)}
			\right)
			-C^{\mathrm{H}_2}
			\sum_{j=1}^{N_{\mathrm{g}}}\sum_{i=1}^{N}F_{i,j,k}
			\Bigg] \Delta t_{\mathrm{upper}} \\
			&+
			\sum_{j=1}^{N_{\mathrm{g}}}\sum_{i=1}^{N}
			\left(
			C^{\mathrm{SU}}b_{i,j,k}^{\mathrm{SU}}
			+
			C^{\mathrm{SD}}b_{i,j,k}^{\mathrm{SD}}
			\right)
			\Bigg\}.
		\end{aligned}
	}
	\label{eq:upper_layer_objective}
\end{equation}
where $C^{\mathrm{W}}$ and $C^{\mathrm{H}2}$ are the electricity and hydrogen prices, respectively; and $C^{\mathrm{SU}}$ and $C^{\mathrm{SD}}$ are the startup and shutdown costs per stack. Power consumption and hydrogen production are integrated over $\Delta t_{\mathrm{upper}}$, whereas startup and shutdown costs are incurred per event.

The upper-layer scheduling problem is summarized as
\begin{equation}
	\begin{aligned}
		\min \quad & J, \\
		\mathrm{s.t.} \quad
		& \eqref{eq:state_exclusive}\text{--}\eqref{eq:start_shutdown_exclusive},\
		\eqref{eq:linear_h2_model}\text{--}\eqref{eq:temperature_power_limit}.
	\end{aligned}
	\label{eq:upper_layer_optimization_summary}
\end{equation}

During rolling operation, the resulting state and power commands are passed to the lower layer, and the executed power and measured thermal states are returned for the next optimization.

\subsection{Performance Assessment Indicators}
\label{subsec:performance_assessment_indicators}

Table~\ref{tab:assessment_indicators} summarizes the indicators used to assess production, thermodynamic performance, startup behavior, degradation, and economics. Hydrogen yield measures total production, energy efficiency measures conversion of electrical input into hydrogen energy, and exergy efficiency further accounts for energy quality. Startup counts and thermal cycling are used to assess degradation, while LCOH measures the resulting life-cycle economic performance.

\begin{table}[tb]
	\scriptsize
	\centering
	\renewcommand{\arraystretch}{1.1}
	\caption{Key assessment indicators for evaluating system operation.}
	\vspace{6pt}
	\label{tab:assessment_indicators}
	\begin{tabularx}{\textwidth}{@{} p{0.20\textwidth} X @{}}
		\hline\hline
		Indicator & Description \\
		\hline
		Hydrogen yield & Total hydrogen production over the investigated horizon, calculated by $V_{\mathrm{H}2,\mathrm{N}}^{\mathrm{tot}}=\int_{0}^{\mathcal{T}}\dot{n}_{\mathrm{H}2}(t)V_{\mathrm{m,N}},\mathrm{d}t$ \cite{xiao2020optimal} \\
		Energy efficiency & Ratio of hydrogen energy output based on the higher heating value (HHV) to net system energy input: $\bar{\eta}_{\mathrm{sys}}=\int_{0}^{\mathcal{T}}\dot{n}_{\mathrm{H}2}(t)\mathit{HHV}_{\mathrm{H}2},\mathrm{d}t/\int_{0}^{\mathcal{T}}P_{\mathrm{sys}}(t),\mathrm{d}t$ \cite{matute2021multi} \\
		Exergy efficiency & Ratio of useful product exergy to total plant input exergy \cite{li2024thermodynamics} \\
		Number of cold startups & Total number of startup events from a low-temperature state after shutdown, denoted by $N_{\mathrm{CS}}=\sum_{y}\sum_{j=1}^{N_{\mathrm{g}}}\sum_{i=1}^{N}N_{\mathrm{CS},i,j,y}$ \cite{ma2025cold} \\
		Number of hot startups & Total number of startup events from a thermally maintained state, denoted by $N_{\mathrm{HS}}=\sum_{y}\sum_{j=1}^{N_{\mathrm{g}}}\sum_{i=1}^{N}N_{\mathrm{HS},i,j,y}$ \cite{ma2025cold} \\
		Voltage degradation & Average annual voltage increase of the stack cluster caused by cold and hot starts \\
		Mechanical damage & Average cumulative fatigue damage of the stack cluster caused by temperature cycling \\
		LCOH & Discounted life-cycle cost divided by hydrogen production corrected for degradation over the project lifetime \cite{xiao2020optimal,matute2021multi,zheng2023offgrid}, as defined in \eqref{eq:lcoh} \\
		\hline\hline
	\end{tabularx}
\end{table}

\subsubsection{Startup Performance}

Cold and hot starts are distinguished by the stack temperature at startup:
\begin{equation}
	\left\{
	\begin{aligned}
		N_{\mathrm{HS},i,j}^{k+1}
		&=
		N_{\mathrm{HS},i,j}^{k}+1
		&& \text{if } T_{i,j}^{k,\mathrm{start}}\geq T_{\mathrm{start}}, \text{ or } 	N_{\mathrm{HS},i,j}^{k}, \text{ otherwise},\\
		N_{\mathrm{CS},i,j}^{k+1}
		&=
		N_{\mathrm{CS},i,j}^{k}+1,
		&& \text{if } T_{i,j}^{k,\mathrm{start}}<T_{\mathrm{start}},  \text{ or } N_{\mathrm{CS},i,j}^{k}, \text{ otherwise}.
	\end{aligned}
	\right.
	\label{eq:start_count}
\end{equation}
where $N_{\mathrm{HS},i,j}^{k}$ and $N_{\mathrm{CS},i,j}^{k}$ are the accumulated hot- and cold-start counts of stack $i$ in group $j$; $T_{i,j}^{k,\mathrm{start}}$ is its temperature at startup; and $T_{\mathrm{start}}$ is the threshold between hot and cold starts.

Startup duration is defined as the time from the startup command to the first instant at which the stack reaches $T_{\mathrm{set}}$:
\begin{equation}
	t_{\mathrm{su}} = t_{\mathrm{set}} - t_{0},
	\label{eq:startup_time}
\end{equation}
where $t_{0}$ is the startup-command time and $t_{\mathrm{set}}$ is the first time at which the stack reaches the thermal-standby target $T_{\mathrm{set}}$.

\subsubsection{Energy and Exergy Efficiencies}

For time $t$ within scheduling interval $k$, the net system power demand includes the actual stack electrochemical input, thermal-standby demand, pump power, and electric-boiler input:
\begin{equation}
		P_{\mathrm{sys}}(t)
		=
		\sum\nolimits_{j=1}^{N_{\mathrm{g}}}\sum\nolimits_{i=1}^{N}
		\left[P_{i,j}^{\mathrm{ele}}(t)+P_{\mathrm{SB}}b_{i,j,k}^{\mathrm{S}}\right]
		+P_{k,\mathrm{pump}}^{(\xi)}
		+P_{k,\mathrm{heat}}^{(\xi)}.
	\label{eq:system_power_demand}
\end{equation}

The instantaneous and average energy efficiencies are
\begin{gather}
	\eta_{\mathrm{sys}}(t)
	=
	{
		\dot{n}_{\mathrm{H}_2}(t)\mathrm{HHV}_{\mathrm{H}_2}
	}/{
		P_{\mathrm{sys}}(t)
	},
	\label{eq:instantaneous_energy_efficiency} \\
	\bar{\eta}_{\mathrm{sys}}(\mathcal{T})
	=
	{\int_{0}^{\mathcal{T}}
		\dot{n}_{\mathrm{H}_2}(t)\mathrm{HHV}_{\mathrm{H}_2}\,\mathrm{d}t
	}\Big/{
		\displaystyle\int_{0}^{\mathcal{T}}
		P_{\mathrm{sys}}(t)\,\mathrm{d}t
	}.
	\label{eq:average_energy_efficiency}
\end{gather}

To account for energy quality, the average exergy efficiency over horizon $\mathcal{T}$ is defined as \cite{li2024thermodynamics}
\begin{equation}
	\bar{\eta}_{\mathrm{ex}}(\mathcal{T})
	={
		\displaystyle\int_{0}^{\mathcal{T}}
		\dot{E}_{x,\mathrm{product}}(t)\,\mathrm{d}t
	} \Big/
	{
		\displaystyle\int_{0}^{\mathcal{T}}
		\dot{E}_{x,\mathrm{in,total}}(t)\,\mathrm{d}t
	}.
	\label{eq:average_exergy_efficiency}
\end{equation}

The total input exergy consists of the electrical exergy supplied to the AWE system, electric boiler, and pumps:
\begin{equation}
	\dot{E}_{x,\mathrm{in,total}}(t)
	=
	\dot{E}_{x,\mathrm{AWE}}(t)
	+
	\dot{E}_{x,\mathrm{boiler}}(t)
	+
	\dot{E}_{x,\mathrm{pump}}(t).
	\label{eq:total_input_exergy}
\end{equation}

The instantaneous exergy efficiency is therefore $\eta_{\mathrm{ex}}(t) = \dot{E}_{x,\mathrm{product}}(t)/\dot{E}_{x,\mathrm{in,total}}(t)$. Detailed material-stream derivations are provided in \ref{sec:appendix_exergy}.

\subsubsection{Degradation Proxies}

Frequent startups and shutdowns cause a small but cumulative increase in cell voltage. For stack $i$ in group $j$, the cell voltage in year $y+1$ is calculated as \cite{li2017sizing,lu2023optimization}
\begin{equation}
{
U_{i,j,y+1}^{\mathrm{cell}}
=
U_{i,j,y}^{\mathrm{cell}}
+
\Delta U_{i,j,y},
}
\label{eq:annual_cell_voltage}
\end{equation}
{where $\Delta U_{i,j,y}$ is the annual startup-shutdown-induced voltage increment of stack $i$ in group $j$, calculated as}
\begin{equation}
{
\Delta U_{i,j,y}
=
N_{\mathrm{CS},i,j,y}\Delta U_{\mathrm{c}}
+
N_{\mathrm{HS},i,j,y}\Delta U_{\mathrm{h}}.
}
\label{eq:voltage_degradation}
\end{equation}
Here, $N_{\mathrm{CS},i,j,y}$ and $N_{\mathrm{HS},i,j,y}$ are the annual cold- and hot-start counts, while $\Delta U_{\mathrm{c}}$ and $\Delta U_{\mathrm{h}}$ are their corresponding voltage increments. The case studies report the mean $\Delta U_{i,j,y}$ over all stacks.

The hydrogen conversion rate in year $y+1$ is calculated as \cite{zheng2023offgrid}
\begin{equation}
	{
		\gamma_{i,j,y+1}
		=
		\left[
		\frac{1}{\gamma_{i,j,y}}
		+
		\Delta U_{i,j,y}
		\frac{2F}{3600M_{\mathrm{H}_2}}
		\right]^{-1},
	}
	\label{eq:hydrogen_conversion_rate}
\end{equation}
where $\gamma{i,j,y}$ is expressed in $\mathrm{kg/kWh}$; $F$ is the Faraday constant; and $M_{\mathrm{H}_2}$ is the hydrogen molar mass in $\mathrm{g/mol}$.

Thermal cycling induced by fluctuating operation causes cumulative mechanical fatigue. For cycle $\ell$ of stack $i$ in group $j$, the temperature range is defined as
\begin{equation}
	\Delta T_{i,j,\ell}
	=
	T_{\max,i,j,\ell}
	-
	T_{\min,i,j,\ell},
	\label{eq:temperature_cycle_range}
\end{equation}
where $T_{\max,i,j,\ell}$ and $T_{\min,i,j,\ell}$ are the maximum and minimum stack temperatures within cycle $\ell$, respectively.

The average cumulative mechanical damage of the stack cluster is quantified as \cite{coffin1954study,zhao2025validity}
\begin{equation}
	\bar{D}
	=
	\frac{1}{N_{\mathrm{g}}N}
	\sum\nolimits_{j=1}^{N_{\mathrm{g}}}
	\sum\nolimits_{i=1}^{N}
	\sum\nolimits_{\ell=1}^{k_{i,j}}
	\left(
	\Delta T_{i,j,\ell} /{\alpha}
	\right)^{\beta},
	\label{eq:mechanical_damage}
\end{equation}
where $k_{i,j}$ is the number of identified temperature cycles and $\alpha$ and $\beta$ are fatigue-model parameters. A smaller $\bar{D}$ indicates less cumulative mechanical damage.

\subsubsection{Economic Assessment}

{Considering cell voltage degradation caused by frequent startups and shutdowns, the levelized cost of hydrogen (LCOH) includes the capital expenditures (CAPEX) of the stack cluster and thermal coupling topology together with annual operating expenditures (OPEX) over the plant lifetime \cite{xiao2020optimal,matute2021multi,superchi2023development,zheng2023offgrid}:}
\begin{equation}
{
\mathrm{LCOH}
=
\frac{
\displaystyle\sum\nolimits_{j=1}^{N_{\mathrm{g}}}\sum\nolimits_{i=1}^{N}\mathrm{CAPEX}_{i,j}
+\mathrm{CAPEX}_{\mathrm{TC}}
+\displaystyle\sum\nolimits_{y=1}^{Y}
\frac{
\displaystyle\sum\nolimits_{j=1}^{N_{\mathrm{g}}}\sum\nolimits_{i=1}^{N}\mathrm{OPEX}_{i,j,y}
+\mathrm{OPEX}_{\mathrm{TC},y}
}{
(1+r)^y
}
}{
\displaystyle\sum\nolimits_{y=1}^{Y}
\sum\nolimits_{j=1}^{N_{\mathrm{g}}}
\sum\nolimits_{i=1}^{N}
\frac{
\gamma_{i,j,y}H_{i,j,y}
}{
\gamma_{i,j,1}(1+r)^y
}
}.
}
\label{eq:lcoh}
\end{equation}


The numerator includes the CAPEX and discounted OPEX of the stack cluster and thermal coupling equipment. Annual OPEX includes electricity for the stacks and thermal auxiliaries together with routine operation and maintenance. $\mathrm{CAPEX}{\mathrm{TC}}$ covers the BHEs, actuators, connecting pipes, and control units, while $\mathrm{OPEX}{\mathrm{TC},y}$ covers their annual operation and maintenance. The denominator is the discounted hydrogen production corrected for startup-shutdown degradation, where $\gamma_{i,j,y}/\gamma_{i,j,1}$ represents the relative hydrogen conversion rate in year $y$, and $H_{i,j,y}$ is the annual hydrogen yield of stack $i$ in group $j$ without degradation.

\section{Results and Discussion}
\label{sec:results_and_discussion}

\subsection{System Parameters}
\label{subsec:system_parameters}

\begin{figure}[tb]
    \centering
    \includegraphics[scale=0.625]{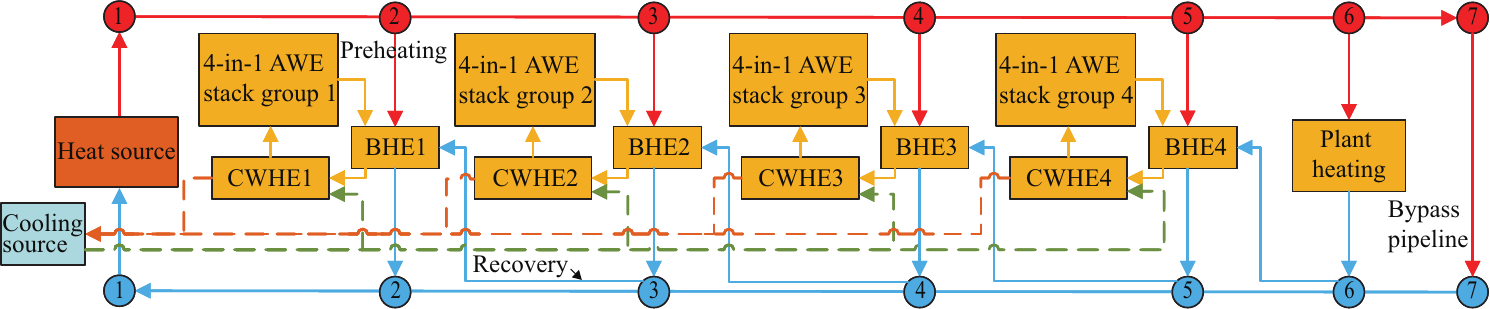}
    \vspace{-6pt}
    \caption{Plant thermal network topology used in the case study.}
    \label{fig:thermal_network_simplified}
\end{figure}

The case study uses design and operational data from an $80~\mathrm{MW}$ hydrogen project in Northern China. The plant contains 16 AWE stacks arranged as four 4-in-1 groups, with each stack rated at $5~\mathrm{MW}$. The AWE parameters are listed in Table~\ref{tab:awe_parameters}, while the plant-building and PTUS parameters are given in Table~\ref{tab:plant_parameters} in \ref{sec:appendix_parameters}. The resulting seven-node thermal network in Fig.~\ref{fig:thermal_network_simplified} contains a PTUS loop and a cooling water loop coupled to the four groups through BHEs and CWHEs, respectively.

The electric boiler acts as the slack heat source at Node~1. The four AWE groups are connected at Nodes~2--5, the plant heating branch at Node~6, and the hydraulic balancing bypass at Node~7. The boiler supply temperature is fixed at $T_{1}^{\mathrm{HS,s}}=363~\mathrm{K}$, while return temperatures are determined from the thermal balance equations. Detailed PTUS pipe parameters are provided in Table~\ref{tab:ptus_parameters} in \ref{sec:appendix_parameters}.

For comparison, the following two methods are considered:

\begin{itemize}
    \item \textbf{Proposed method (PM):} PM applies the framework in Section~\ref{sec:two_layer_control}, coordinating waste heat recovery, preheating, and thermal standby across the plant. The upper layer uses a $24~\mathrm{h}$ rolling horizon with a $15~\mathrm{min}$ step, while the lower layer regulates current and temperature every $15~\mathrm{s}$.

    \item \textbf{Benchmark method (BM):} BM represents conventional independent thermal management. The scheduling layer uses the same renewable power and stack operating constraints as PM, while the lower layer uses conventional lye circulation before startup and regulates the four 4-in-1 groups independently through lye flow and active cooling.
\end{itemize}

Measured wind power is used for the weekly assessment, while photovoltaic (PV) data are used for the seasonal and annual studies. Both methods use identical renewable-power inputs, plant capacities, and common equipment limits. The models are implemented in \textit{Wolfram Mathematica 13.0}, and the MILP problems are solved using \textit{Gurobi 11.0.3}.

\subsection{Dynamic Performance Analysis}
\label{subsec:dynamic_performance_thermal_strategy}

Dynamic performance is first examined through single-stack startup and shutdown, thermal standby in $N$-in-1 groups, and plant air temperature regulation.

\subsubsection{Analysis of Startup and Shutdown of a Single Stack}
\label{subsubsec:single_stack_analysis}

\begin{figure}[tb]
	\centering
	\includegraphics[scale=1]{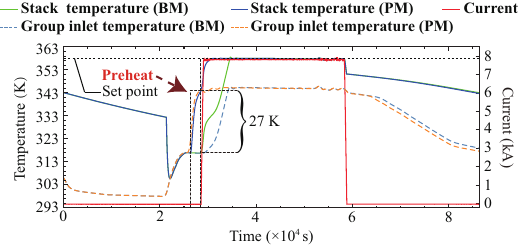}
	\vspace{-6pt}
	\caption{Stack current and temperature responses during single-stack startup and shutdown under PM and BM.}
	\label{fig:single_stack_curves}
\end{figure}

Fig.~\ref{fig:single_stack_curves} compares stack current and temperature during startup and shutdown under PM and BM. Under BM, the lye flow is increased $2.5~\mathrm{h}$ before startup. The initial circulation removes heat and lowers the stack temperature, after which the warmer return lye gradually raises the group inlet temperature from $20{,}500$ to $28{,}000~\mathrm{s}$. During startup, electrochemical heat further increases the inlet and stack temperatures. Active cooling begins when the stack temperature exceeds $348~\mathrm{K}$, and the temperature approaches $358~\mathrm{K}$ after about $21~\mathrm{min}$. During shutdown, lye and cooling-water circulation continue briefly to remove residual gas and heat before natural cooling.

This slow thermal buildup restricts startup under BM. Under PM, the BHE supplies PTUS heat to the lye loop beginning $2{,}400~\mathrm{s}$ before startup. The stack therefore starts up at a temperature $27~\mathrm{K}$ higher than under BM, reducing the startup duration from $7{,}000$ to $2{,}500~\mathrm{s}$, or by $64.3\%$.

\subsubsection{Thermal Standby Process in Stack Clusters}
\label{subsubsec:thermal_standby_analysis}

\begin{figure}[tb]
	\centering
	\includegraphics[width=0.85\textwidth]{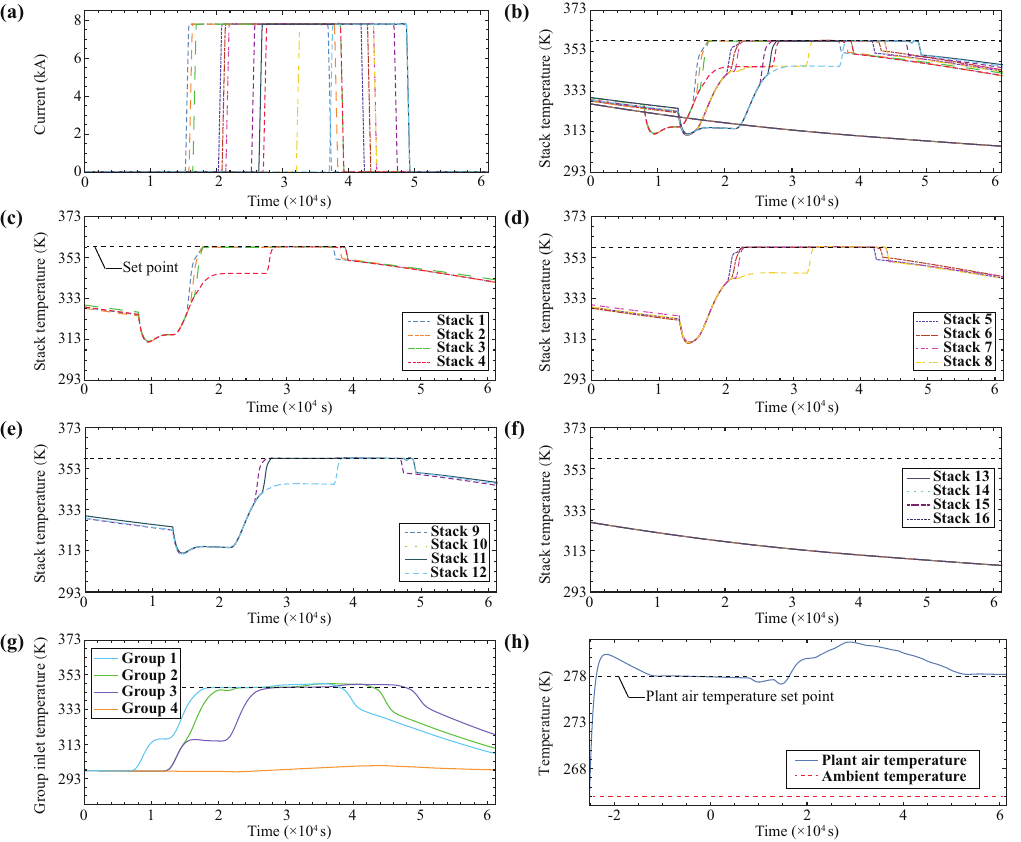}\vspace{-6pt}
	\caption{Temperature dynamics of the hydrogen plant during a daily production cycle under PM. (a) Currents of 16 stacks; (b) Stack temperatures of all groups; (c) Stack temperatures of group~1; (d) Stack temperatures of group~2; (e) Stack temperatures of group~3; (f) Stack temperatures of group~4; (g) Group inlet temperatures of groups~1--4; (h) Plant air temperature.}
	\label{fig:group_response}
\end{figure}

Fig.~\ref{fig:group_response} shows the current and temperature responses of the four 4-in-1 groups under PM. The current profiles represent typical startup and shutdown sequences and are used to examine intra- and inter-group thermal standby.

As shown in Figs.~\ref{fig:group_response}(c), (d), and (g), groups~1 and~2 undergo intra-group thermal standby during $16{,}500$--$22{,}000~\mathrm{s}$. Before this period, lye circulation raises the group inlet temperature and narrows the temperature differences among the stacks. After about $16{,}500~\mathrm{s}$, waste heat from high-load production stacks is redistributed through the shared lye loop to maintain the other stacks in thermal standby. After $28{,}000~\mathrm{s}$, three stacks in each group operate at rated power, while the fourth remains near the standby temperature through intra-group heat redistribution.

For group~3, Figs.~\ref{fig:group_response}(e) and (g) show that no stack operates at high load before startup, so local heat redistribution is insufficient. Inter-group thermal standby is therefore established through the PTUS by about $22{,}000~\mathrm{s}$. With the stack operating states and power levels temporarily unchanged, the subsequent heat transfer raises the group inlet temperature from about $313~\mathrm{K}$ to $348~\mathrm{K}$ by $27{,}000~\mathrm{s}$. Group~4 does not enter thermal standby, and its stack temperatures gradually approach the plant air temperature, as shown in Fig.~\ref{fig:group_response}(f).

\subsubsection{Plant Air Temperature Dynamics and Regulation}
\label{subsubsec:plant_temperature_dynamics}

Fig.~\ref{fig:group_response}(h) shows the plant air temperature. Before $t=0$, the PTUS preheats the building through radiators to approximately the $278~\mathrm{K}$ set point. After stack~1 enters high-load operation at about $16{,}500~\mathrm{s}$, equipment heat dissipation raises the plant air temperature even with the radiator valve at its minimum opening. The temperature peaks near $30{,}000~\mathrm{s}$ and then decreases as heat is removed through the building envelope and ventilation, returning to approximately $278~\mathrm{K}$ by $55{,}000$--$60{,}000~\mathrm{s}$.

\subsection{Component Energy and Exergy Analysis}
\label{subsec:sectional_energy_exergy_analysis}

\begin{figure}[!htbp]
	\centering
	\includegraphics[width=\textwidth]{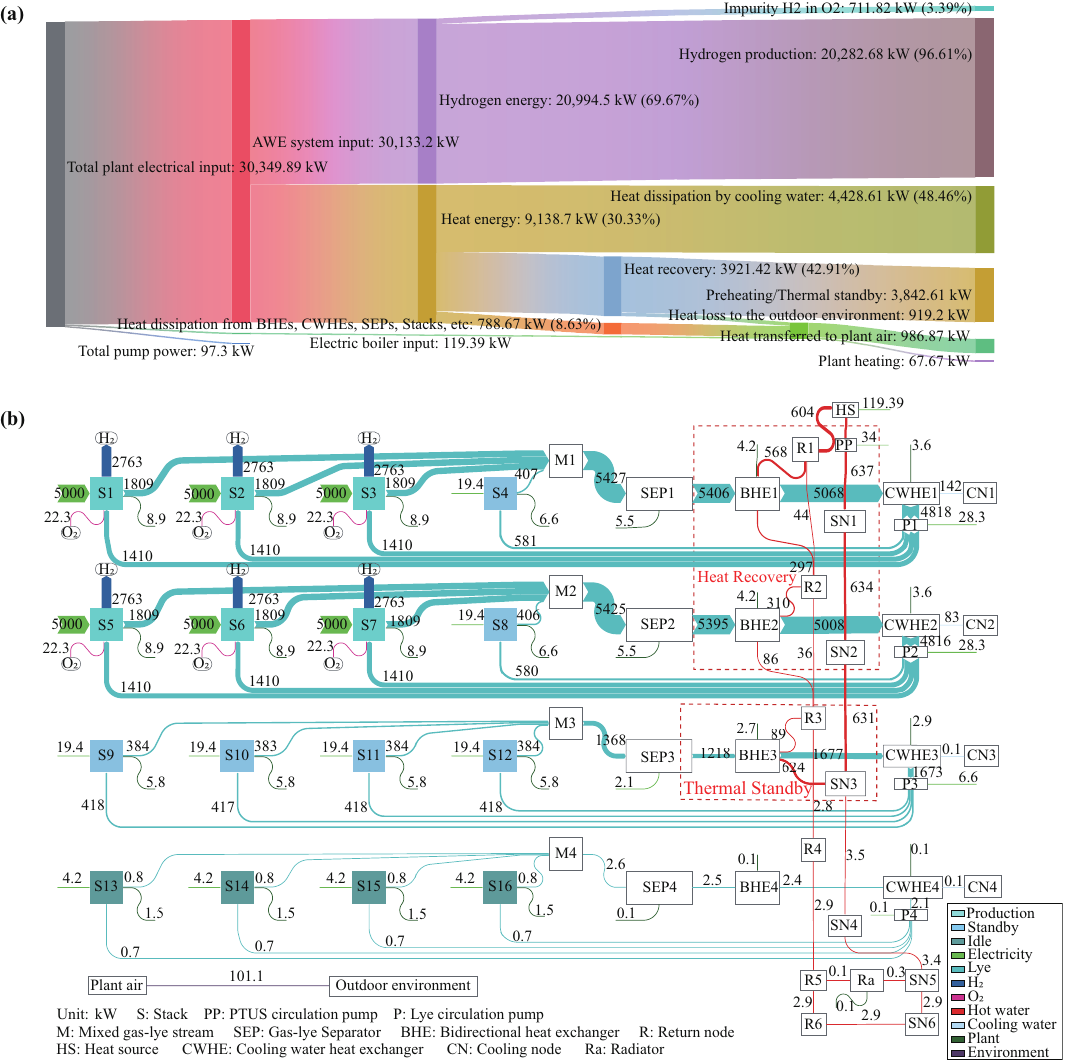}
	\caption{Instantaneous energy and exergy flow diagrams of the hydrogen plant at $t=22{,}000~\mathrm{s}$. (a) Energy flow diagram of the hydrogen plant; (b) Sankey diagram of process exergy flow.}
	\label{fig:energy_exergy_sankey}
\end{figure}

The operating snapshot at $t=22{,}000~\mathrm{s}$ in Section~\ref{subsubsec:thermal_standby_analysis} is selected for instantaneous energy and exergy analysis. At this time, groups~1 and~2 each have three stacks in production and one in intra-group thermal standby, group~3 is under inter-group thermal standby, and group~4 is idle. Because the stack operating states and electrical powers are temporarily unchanged, this instant is treated as a quasi-steady snapshot for component-level energy and exergy accounting.

\subsubsection{Energy Flow Analysis}
\label{subsubsec:energy_consumption_analysis}

At $t=22{,}000~\mathrm{s}$, the total plant electrical input is $30{,}349.89~\mathrm{kW}$, including $30{,}133.2~\mathrm{kW}$ for the AWE system, $97.3~\mathrm{kW}$ for pumps, and $119.39~\mathrm{kW}$ for the electric boiler, as shown in Fig.~\ref{fig:energy_exergy_sankey}. The boiler produces $117.00~\mathrm{kW}$ of heat. Coordinated power and lye flow regulation keeps the stacks in production within a favorable electrochemical range.

Of the AWE electrical input, $20{,}994.5~\mathrm{kW}$ is converted into hydrogen chemical energy, giving an AWE electrical-to-hydrogen efficiency of $69.67\%$. The remaining $9{,}138.7~\mathrm{kW}$ appears as heat in the gas-lye system. Of this heat, $4{,}428.61~\mathrm{kW}$ ($48.46\%$) is removed by cooling water, $3{,}921.42~\mathrm{kW}$ ($42.91\%$) is recovered, and $788.67~\mathrm{kW}$ ($8.63\%$) is dissipated from the stacks, SEPs, CWHEs, BHEs, and other BoP equipment. Part of this dissipated heat contributes to plant heating, while excess heat is rejected through the building envelope and ventilation.

Of the $3{,}921.42~\mathrm{kW}$ recovered heat, $3{,}842.61~\mathrm{kW}$ is used for stack preheating and thermal standby, while $78.81~\mathrm{kW}$ supplies plant heating and PTUS heat losses. Because the recovered heat alone cannot maintain the required PTUS supply temperature of $363~\mathrm{K}$, the electric boiler contributes another $117.00~\mathrm{kW}$. Together, recovered heat and boiler output provide $4{,}038.42~\mathrm{kW}$. If this thermal demand were supplied entirely by the electric boiler, its electrical input would be $4{,}120.84~\mathrm{kW}$. The proposed thermal coupling therefore reduces the equivalent electric-boiler input by $4{,}001.45~\mathrm{kW}$ at this operating point by reusing low-grade stack heat.

\subsubsection{Exergy Analysis}
\label{subsubsec:results_exergy_analysis}

\begin{table}[tb]\scriptsize\centering
    \renewcommand{\arraystretch}{1.05}
    \caption{Instantaneous exergy destruction of the main plant processes.}
    \vspace{6pt}
    \label{tab:exergy_destruction}
    \begin{tabular}{cccccc}
        \hline\hline
        Process & \tabincell{c}{Exergy destruction ($\mathrm{kW}$)}
        & \tabincell{c}{Proportion ($\%$)}
        & Process & \tabincell{c}{Exergy destruction ($\mathrm{kW}$)}
        & \tabincell{c}{Proportion ($\%$)} \\
        \hline
        Stacks & $11{,}415.2$ & ${82.89}$
        & BHE & $294.9$ & $2.14$ \\

        Electric boiler & ${72.59}$ & ${0.53}$
        & CWHE & $210.9$ & $1.53$ \\

        Lye mixing & $982.6$ & ${7.13}$
        & Cooling water & $225.2$ & $1.64$ \\

        SEP & $187.9$ & $1.36$
        & Total pumps & ${97.3}$ & $0.71$ \\

        Others & $285.2$ & $2.07$
        & \tabincell{c}{Total exergy\\destruction}
        & ${13{,}771.79}$ & $100.00$ \\
        \hline\hline
    \end{tabular}
\end{table}

At $t=22{,}000~\mathrm{s}$, the total exergy input is $30{,}349.89~\mathrm{kW}$ and the plant exergy efficiency is $54.62\%$. Table~\ref{tab:exergy_destruction} summarizes the main sources of exergy destruction. As shown in Fig.~\ref{fig:energy_exergy_sankey}(b), AWE stacks dominate, with $11{,}415.2~\mathrm{kW}$ of electrochemical irreversibility, or $82.89\%$ of the total. Lye mixing and separation contribute $982.6~\mathrm{kW}$ ($7.13\%$) and $187.9~\mathrm{kW}$ ($1.36\%$), respectively. Heat transfer across finite temperature differences causes $210.9~\mathrm{kW}$ ($1.53\%$) of exergy destruction in the CWHEs and $294.9~\mathrm{kW}$ ($2.14\%$) in the BHEs. The electric boiler contributes only $72.59~\mathrm{kW}$ ($0.53\%$) because most thermal demand is met by recovered stack heat.

The thermal utilization pathways provide $830~\mathrm{kW}$ of low-grade thermal exergy, including $482~\mathrm{kW}$ through PTUS-mediated waste heat recovery and $348~\mathrm{kW}$ through intra-group redistribution. These flows support the thermal standby of six stacks and reduce the use of high-grade electrical exergy. The exergy benefit of PM therefore arises mainly from recovering and redistributing stack waste heat that would otherwise be rejected.

\subsection{Weekly Energy Management under Wind-Power Fluctuations}
\label{subsec:weekly_wind_energy_management}

\begin{figure}[tb]
  \centering
\includegraphics[scale=0.9]{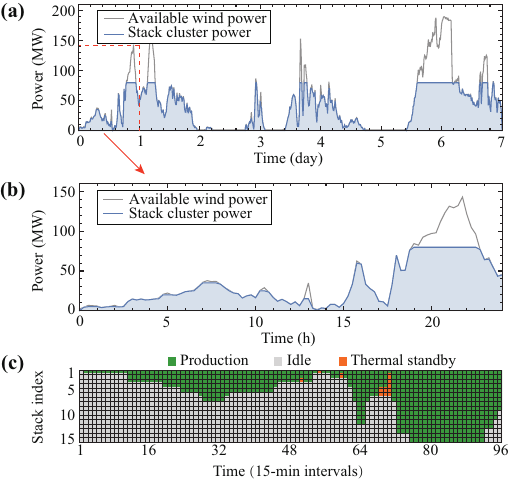}
\vspace{-6pt}
\caption{Wind power operation over seven days. (a) Wind power availability and stack-cluster power; (b) First-day power profiles; and (c) First-day stack operating states.}
  \label{fig:7_day_schedule}
\end{figure}

\begin{figure}[htbp]
  \centering
\includegraphics[width=0.95\textwidth,height=0.88\textheight,keepaspectratio]{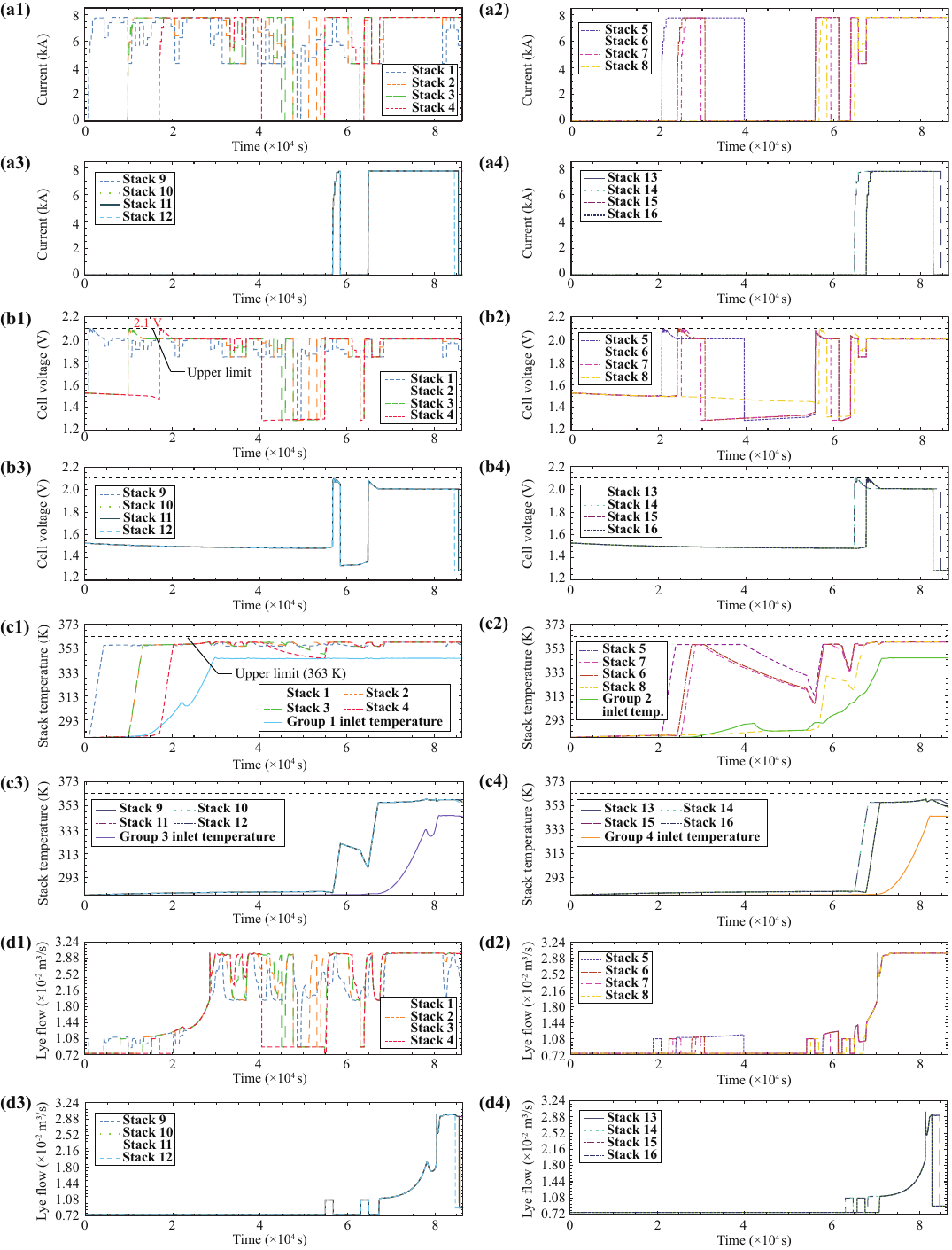}
\vspace{-6pt}
  \caption{Thermal management process under BM. (a1--a4) Stack current; (b1--b4) Cell voltage; (c1--c4) Group inlet and stack temperatures; (d1--d4) Lye flow rates.}
  \label{fig:weekly_bm_overview}
\end{figure}

\begin{figure}[htbp]
  \centering
\includegraphics[width=0.95\textwidth,height=0.88\textheight,keepaspectratio]{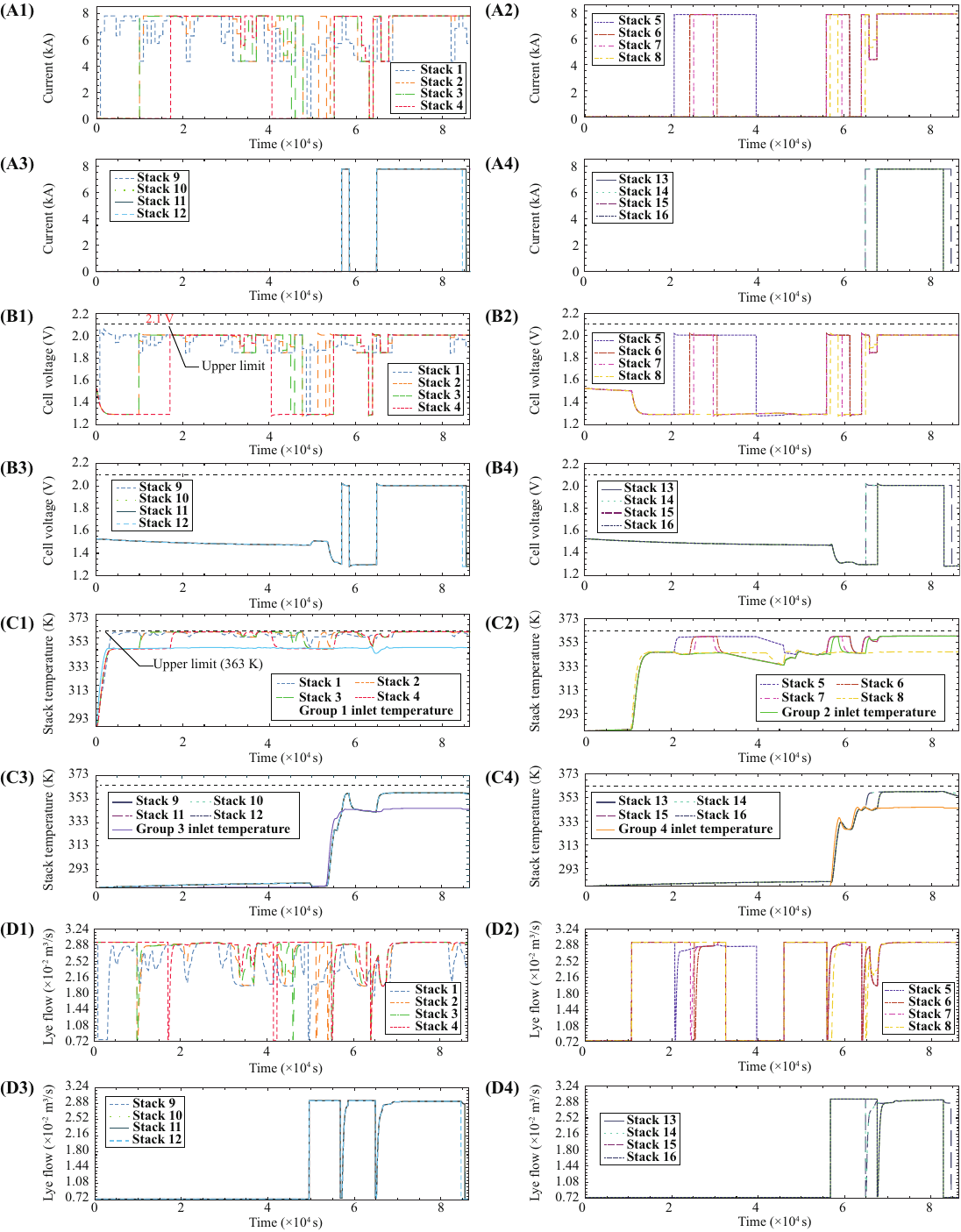}
\vspace{-6pt}
  \caption{Thermal management process under PM. (A1--A4) Stack current; (B1--B4) Cell voltage; (C1--C4) Group inlet and stack temperatures; (D1--D4) Lye flow rates.}
  \label{fig:weekly_pm_overview}
\end{figure}

A continuous seven-day simulation using measured wind power evaluates operation under renewable fluctuations. The upper layer updates unit commitment and stack power commands every $15~\mathrm{min}$, while the lower layer performs $15~\mathrm{s}$ feedback regulation. The first day is examined in detail. Based on the power and operating-state profiles in Figs.~\ref{fig:7_day_schedule}(b) and (c), Figs.~\ref{fig:weekly_bm_overview} and~\ref{fig:weekly_pm_overview} compare stack currents, cell voltages, temperatures, and lye flows under BM and PM.

\begin{figure}[tb]
  \centering
\includegraphics[width=0.85\textwidth]{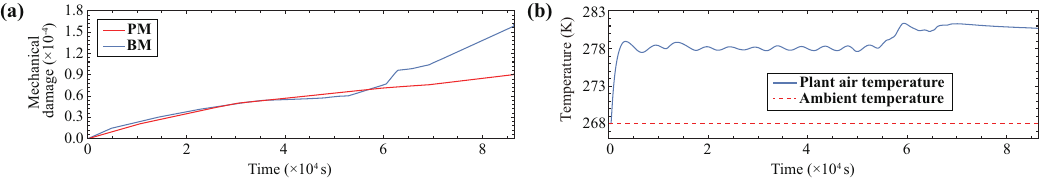}
\vspace{-6pt}
  \caption{First-day simulation results. (a) Average cumulative mechanical damage of the 16 stacks under PM and BM; (b) Plant air and ambient temperature variations under PM.}
  \label{fig:degradation_curves}
\end{figure}

\subsubsection{Preheating, Startup, and Current Increase}
\label{subsubsec:weekly_preheating_startup}

During $0$--$20{,}000~\mathrm{s}$, the upper layer commands group 1 to start, while stacks~2--4 are placed in intra-group thermal standby under PM. The lye pumps initially operate at $v_{\mathrm{low}}=0.25v_{\mathrm{lye,rate}}$ ($0.72 \times 10^{-2}~\mathrm{m}^3/\mathrm{s}$) and then follow the preheating command. PTUS heat raises the group inlet temperature from $278~\mathrm{K}$ to $323~\mathrm{K}$ before startup and to about $344~\mathrm{K}$ by $t=2{,}500~\mathrm{s}$, as shown in Fig.~\ref{fig:weekly_pm_overview}(C1) and (D1).

During $900$--$20{,}000~\mathrm{s}$, stacks~1--4 successively enter production. Their higher startup temperatures reduce cell voltage and allow the stack currents to increase more rapidly within the voltage limit. Stacks~2--4 consequently reach the rated current of $7{,}800~\mathrm{A}$ within about 2 min, as shown in Fig.~\ref{fig:weekly_pm_overview}(A1) and (B1).

Under BM, the group inlet temperature remains below $313~\mathrm{K}$ for longer because pre-startup lye circulation provides no external heat input. The resulting higher cell voltage imposes a tighter temperature-dependent current limit, as shown in Fig.~\ref{fig:weekly_bm_overview}(a1)--(d1). After the stacks reach normal operation, the temperature differences between the two methods decrease. However, the first-day cumulative mechanical damage reaches approximately $0.8 \times 10^{-4}$ under PM and $1.54 \times 10^{-4}$ under BM, as shown in Fig.~\ref{fig:degradation_curves}(a).

\subsubsection{Operation at High Load}
\label{subsubsec:weekly_high_load_operation}

During $22{,}000$--$40{,}000~\mathrm{s}$, stacks~5--7 start up as renewable power increases, as shown in Fig.~\ref{fig:weekly_pm_overview}(A2). Waste heat from high-load stacks~1--4 is transferred through the PTUS to preheat group~2, raising its inlet temperature to about $345~\mathrm{K}$. With the lye flow increased to $v_{\mathrm{pre}}=v_{\mathrm{lye,rate}}$ during preheating and then regulated during production, stacks~5--8 approach the hot-startup condition, as shown in Figs.~\ref{fig:weekly_pm_overview}(C2) and (D2). Waste heat from stacks~5--7 subsequently maintains stack~8 in thermal standby.

Inter-group thermal standby is also activated for stacks~5--8 around $49{,}000~\mathrm{s}$, stacks~9--11 around $57{,}000~\mathrm{s}$, and stacks~13--16 during $57{,}000$--$65{,}000~\mathrm{s}$. Intra-group standby occurs for stacks~8 and~16. Together, these heat-sharing paths slow cooling and reduce temperature cycling relative to BM.

\subsubsection{Current Reduction and Shutdown}
\label{subsubsec:weekly_current_reduction_shutdown}

As renewable power declines, stack currents decrease under both methods. Under PM, partly idle groups remain thermally supported through intra-group circulation or the PTUS. For example, when group~3 is idle during $59{,}000$--$65{,}000~\mathrm{s}$, waste heat from other high-load groups slows its temperature decrease, as shown in Fig.~\ref{fig:weekly_pm_overview}(C3).

Under BM, the absence of heat redistribution leads to faster cooling after shutdown. The stacks in group~3 reach only about $356~\mathrm{K}$ during high-load operation and cool rapidly once the lye pump is stopped. Restarting circulation before the next startup introduces colder lye and causes a further temperature drop, increasing the subsequent cold-startup burden.

The plant air temperature follows the mechanism discussed in Section~\ref{subsubsec:plant_temperature_dynamics}. The short rise after $55{,}000~\mathrm{s}$ results from the startup of groups~3 and~4 and is subsequently limited by increased ventilation.

\subsubsection{Weekly Performance Evaluation}
\label{subsubsec:weekly_evaluation}

Table~\ref{tab:weekly_performance_comparison} summarizes the performance over the continuous seven-day simulation.

\begin{table}[tb]\scriptsize\centering
	\renewcommand{\arraystretch}{1.1}
	\caption{Weekly performance comparison under PM and BM.}\vspace{6pt}
	\label{tab:weekly_performance_comparison}
	\begin{tabular}{cccccccc}
		\hline\hline
		Method &
		\tabincell{c}{$\mathrm{H}_2$ yield\\($\mathrm{kNm}^3$)} &
		\tabincell{c}{Energy\\efficiency ($\%$)} &
		\tabincell{c}{Exergy\\ efficiency ($\%$)} &
		\tabincell{c}{LCOH\\($\mathrm{CNY/kg}$)} &
		\tabincell{c}{Cold/hot\\startups} &
		\tabincell{c}{Mech.\ damage\\($\times 10^{-2}$)} &
		\tabincell{c}{Annualized voltage\\degradation ($\mathrm{V}$)} \\
		\hline
		PM & $1{,}025.1$ & $66.695$ & $56.43$ & $33.391$ & $0$/$206$ & $0.076$ & $0.005356$ \\
		BM & $1{,}008.7$ & $65.763$ & $52.15$ & $34.038$ & $73$/$133$ & $0.120$ & $0.022438$ \\
		\hline\hline
	\end{tabular}

\vspace{6pt}
\parbox{0.98\textwidth}{\footnotesize\emph{Note:} {
				The weekly cold- and hot-startup counts are annualized by a factor of $52$ before calculating $\Delta U_{i,j,y}$ from \eqref{eq:voltage_degradation}; the reported voltage degradation is averaged over all stacks.}}
	
\end{table}

\emph{a) Production and efficiency:} PM increases weekly hydrogen yield by $1.63\%$. Energy and exergy efficiencies increase by $0.932$ and $4.28$ percentage points, respectively. These gains arise from shorter startup delays, higher allowable current at elevated stack temperatures, and reduced electric heating through waste heat recovery.

\emph{b) Economics:} PM reduces LCOH from $34.038$ to $33.391~\mathrm{CNY/kg}$, or by $1.90\%$. The reduction results from higher hydrogen production, lower simulated electricity consumption, and smaller startup-shutdown-induced conversion losses.

\emph{c) Startup and degradation:} PM eliminates the $73$ cold startups observed under BM, converting all startup events to hot startups. Average cumulative mechanical damage decreases by $36.67\%$, while annualized voltage degradation decreases by $76.13\%$. These reductions result from shallower temperature cycling and higher stack temperatures before startup.

\subsection{Seasonal Scenario Analysis and Annual Performance}
\label{subsec:seasonal_scenario_analysis}

Long-term performance is evaluated over a representative 364-day year consisting of four 91-day seasons. Fig.~\ref{fig:plant_temp_curve} and Table~\ref{tab:seasonal_annual_comparison} summarize the seasonal and annual results.

\subsubsection{Seasonal Performance}

The benefit of PM varies with renewable power and ambient temperature. Spring produces the most hydrogen because of favorable renewable availability and temperature conditions. In summer, the smaller cold-startup penalty reduces the relative benefit of preheating. Autumn has lower hydrogen production because of reduced renewable availability, although relatively stable operation maintains high energy efficiency. Winter has the highest heating and startup demand and therefore shows the largest gain in exergy efficiency.

\begin{figure}[!tb]
	\centering
	\includegraphics[scale=0.85]{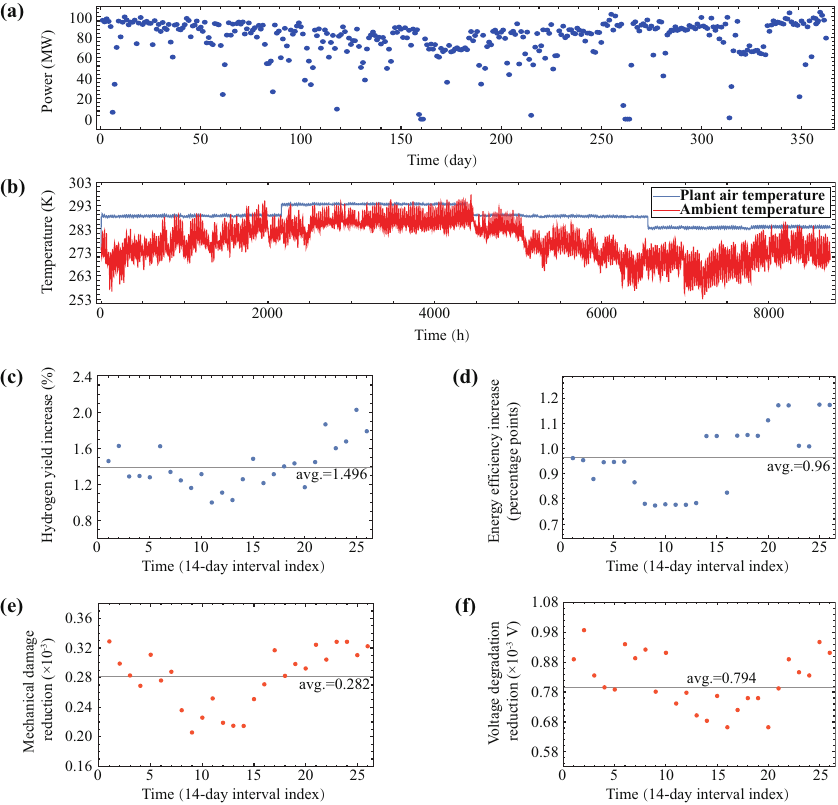}
	\caption{Annual operating performance. (a) Daily peak PV power; (b) Plant air and ambient temperatures under PM;  (c) 14-day hydrogen yield improvements of PM over BM; (d) Energy efficiency increase; (e) Mechanical damage reduction; and (f) Voltage degradation reduction.}
	\label{fig:plant_temp_curve}
\end{figure}

Across the full year, PM increases hydrogen yield by $1.50\%$ and raises energy and exergy efficiencies by $0.99$ and $4.33$ percentage points, respectively. The averages in Fig.~\ref{fig:plant_temp_curve}(c)--(f) are unweighted means of the 26 individual 14-day intervals and therefore differ slightly from the annual total-based indicators in Table~\ref{tab:seasonal_annual_comparison}. Overall, PM retains its performance advantage across all seasons despite changes in renewable availability and ambient temperature.

\begin{table}[tb]\scriptsize\centering
\renewcommand{\arraystretch}{1.14}
\setlength{\tabcolsep}{3pt}
\caption{Seasonal and annual operating performance under PM and BM.}
\vspace{6pt}
\label{tab:seasonal_annual_comparison}
\begin{adjustbox}{max width=\textwidth}
	\begin{tabular}{ccccccccc}
		\hline\hline
		Season
		& Method
		& \tabincell{c}{$\mathrm{H}_2$ yield\\($\mathrm{kNm}^3$)}
		& \tabincell{c}{Energy \\efficiency ($\%$)}
		& \tabincell{c}{Exergy \\ efficiency ($\%$)}
		& \tabincell{c}{LCOH\\($\mathrm{CNY/kg}$)}
		& \tabincell{c}{Cold/hot\\startups}
		& \tabincell{c}{Mechanical \\ damage ($\times 10^{-2}$)}
		& \tabincell{c}{Annualized voltage\\degradation ($\mathrm{V}$)} \\
		\hline
		
		\multirow{2}{*}{Spring}
		& PM & $7{,}249.27$ & $66.49$ & $52.31$ & $33.471$ & $0/1{,}249$ & $0.445$ & $0.002498$ \\
		& BM & $7{,}145.47$ & $65.53$ & $48.29$ & $34.656$ & $1{,}249/0$ & $0.629$ & $0.024980$ \\
		\hline
		
		\multirow{2}{*}{Summer}
		& PM & $7{,}187.33$ & $66.56$ & $53.63$ & $33.454$ & $0/1{,}270$ & $0.384$ & $0.002540$ \\
		& BM & $7{,}103.31$ & $65.78$ & $50.54$ & $34.444$ & $1{,}146/124$ & $0.541$ & $0.023168$ \\
		\hline
		
		\multirow{2}{*}{Autumn}
		& PM & $6{,}225.05$ & $66.57$ & $53.51$ & $33.537$ & $0/1{,}020$ & $0.299$ & $0.002040$ \\
		& BM & $6{,}126.76$ & $65.52$ & $49.22$ & $34.853$ & $985/35$ & $0.486$ & $0.019770$ \\
		\hline
		
		\multirow{2}{*}{Winter}
		& PM & $6{,}915.57$ & $66.55$ & $52.85$ & $33.741$ & $0/1{,}253$ & $0.387$ & $0.002506$ \\
		& BM & $6{,}793.47$ & $65.38$ & $46.78$ & $34.694$ & $1{,}208/45$ & $0.592$ & $0.024250$ \\
		\hline
		
		\multirow{2}{*}{Annual}
		& PM & $27{,}577.22$ & $66.54$ & $52.96$ & $33.541$ & $0/4{,}792$ & $1.518$ & $0.002396$ \\
		& BM & $27{,}169.01$ & $65.55$ & $48.63$ & $34.657$ & $4{,}588/204$ & $2.251$ & $0.023042$ \\
		\hline\hline
	\end{tabular}
\end{adjustbox}

\vspace{6pt}
\parbox{0.98\textwidth}{\footnotesize\emph{Note:} {For the seasonal cases, the cold- and hot-startup counts are annualized by a factor of $4$ before calculating $\Delta U_{i,j,y}$ from \eqref{eq:voltage_degradation}; the reported voltage degradation is averaged over all stacks, whereas the annual value is calculated directly from the annual startup counts before averaging. Annual mechanical damage is recalculated from the continuous annual temperature trajectories using \eqref{eq:mechanical_damage}. Annual energy and exergy efficiencies are evaluated from annual totals, and LCOH is calculated from \eqref{eq:lcoh}.}}
\end{table}

\subsubsection{Annual Degradation and Economics}

The annual simulation further demonstrates the long-term benefits of PM. Results in Table~\ref{tab:seasonal_annual_comparison} show that BM produces $4{,}588$ cold startups and $204$ hot startups, whereas PM converts all $4{,}792$ startup events to hot startups. Average cumulative mechanical damage consequently decreases from $2.251\times10^{-2}$ to $1.518\times10^{-2}$, a reduction of $32.56\%$. Annual voltage degradation decreases from $0.023042$ to $0.002396~\mathrm{V}$, or by $89.60\%$. PTUS preheating and thermal standby maintain higher stack temperatures before startup and reduce deep temperature cycling throughout the year.

Together with the $1.50\%$ increase in annual hydrogen yield, the reduction in electricity consumption and startup-shutdown-induced conversion loss lowers LCOH from $34.657$ to $33.541~\mathrm{CNY/kg}$, a reduction of $3.22\%$. The annual results therefore show that plant-wide thermal coordination improves production, thermodynamic efficiency, degradation performance, and economics, with the strongest thermal benefit under cold operating conditions.

\section{Conclusion}
\label{sec:conclusion}

This paper develops a plant-wide thermal model and hierarchical electricity-heat coordination framework for large-scale cold-region ReP2H plants. The framework couples the AWE stack cluster, BoP, PTUS, and plant building through bidirectional heat exchange, while coordinating minute-scale production scheduling with second-scale current and thermal regulation. Simulations of an $80~\mathrm{MW}$ ReP2H plant lead to the following findings:

\begin{enumerate}
	\item Low stack temperature prolongs startup and tightens the current limit imposed by the cell-voltage constraint. Thermal standby is effective for short shutdowns, while PTUS-assisted preheating is needed before low-temperature restart. In the single-stack case, the proposed method reduces startup duration by $64.3\%$.
	
	\item In $N$-in-1 AWE groups, waste heat can be redistributed within a group or transferred through the PTUS among groups for thermal standby and preheating, while surplus heat can support plant heating before rejection. The heat-transfer mode is determined by stack thermal states, scheduled startup times, available waste heat, and PTUS conditions.
	
	\item The benefits of thermal coordination accumulate from individual startup events to annual operation. In the evaluated annual scenario, the proposed method increases hydrogen yield by $1.50\%$ and raises energy and exergy efficiencies by $0.99$ and $4.33$ percentage points, respectively. Mechanical damage, startup-shutdown-induced voltage degradation, and LCOH decrease by $32.56\%$, $89.60\%$, and $3.22\%$, respectively.
\end{enumerate}

Future work will validate the proposed framework in pilot-scale or industrial operation, improve electrochemical and mechanical degradation models, and incorporate renewable-power and heat-demand uncertainty into scheduling. The framework can also be extended to hydrogen storage and downstream processes for broader plant-wide energy coordination.

\appendix

\section{Seasonal Fitting of Auxiliary Power for Thermal Management}
\label{sec:appendix_fitting}

\setcounter{figure}{0}
\renewcommand{\thefigure}{A\arabic{figure}}
\setcounter{table}{0}
\renewcommand{\thetable}{A\arabic{table}}

\begin{figure}[H]
    \centering
\includegraphics[scale=0.9]{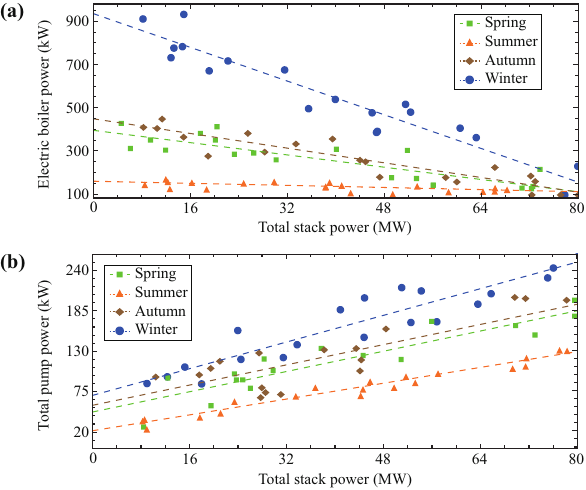}
\caption{Simulated and fitted thermal management auxiliary power under different seasons. (a) Electric boiler input power; (b) Total pump power.}
    \label{fig:appendix_power_fit}
\end{figure}

Fig.~\ref{fig:appendix_power_fit} compares the simulated and fitted electric-boiler and total pump powers against total stack power. The seasonal fitting coefficients obtained from the detailed thermal simulations are listed in Table~\ref{tab:seasonal_fitting_params}.

\begin{table}[H]\scriptsize\centering
    \renewcommand{\arraystretch}{1.05}
    \caption{Seasonal fitting parameters of thermal management auxiliary power.}\vspace{6pt}
    \label{tab:seasonal_fitting_params}
    \begin{tabular}{ccccc}
        \hline\hline
        Season & \tabincell{c}{{$\alpha_{\mathrm{pump}}^{(\xi)}$}} & \tabincell{c}{{$b_{\mathrm{pump}}^{(\xi)}$} ($\mathrm{kW}$)} & \tabincell{c}{{$\alpha_{\mathrm{heat}}^{(\xi)}$}} & \tabincell{c}{{$b_{\mathrm{heat}}^{(\xi)}$} ($\mathrm{kW}$)} \\
        \hline
        Spring & $1.723 \times 10^{-3}$ & $47.3$ & $-3.558 \times 10^{-3}$ & $395.5$ \\
        Summer & $1.347 \times 10^{-3}$ & $21.8$ & $-6.048 \times 10^{-4}$ & $159.8$ \\
        Autumn & $1.720 \times 10^{-3}$ & $56.4$ & $-4.249 \times 10^{-3}$ & $448.9$ \\
        Winter & $2.269 \times 10^{-3}$ & $69.9$ & $-9.758 \times 10^{-3}$ & $936.0$ \\
        \hline\hline
    \end{tabular}
\end{table}

\section{Exergy Evaluation}
\label{sec:appendix_exergy}

The reference environment is defined by $T_0=298.15~\mathrm{K}$ and $p_0=1~\mathrm{atm}$. Kinetic and potential exergy are neglected because they are small relative to thermal and chemical exergy. The total exergy rate of a material stream is therefore the sum of its physical and chemical exergy rates \cite{li2024thermodynamics}:
\begin{equation}
	\dot{E}_{x,m}
	=
	\dot{E}_{x,m}^{\mathrm{ph}}
	+
	\dot{E}_{x,m}^{\mathrm{ch}}.
	\label{eq:appendix_total_material_exergy}
\end{equation}

The physical exergy rate associated with departure from the reference state is
\begin{equation}
	\dot{E}_{x,m}^{\mathrm{ph}}
	=
	\dot{m}
	\left[
	\left(h-h_0\right)
	-
	T_0\left(s-s_0\right)
	\right],
	\label{eq:appendix_physical_exergy}
\end{equation}
where $\dot{m}$ is the mass flow rate; $h$ and $s$ are the specific enthalpy and entropy of the stream; and $h_0$ and $s_0$ are their values at the reference state.

The chemical exergy rate associated with stream composition is expressed as
\begin{equation}
	\dot{E}_{x,m}^{\mathrm{ch}}
	=
	\sum\nolimits_q
	\dot{n}_q
	e_{x,q}^{\mathrm{ch}},
	\label{eq:appendix_chemical_exergy}
\end{equation}
where $\dot{n}q$ and $e{x,q}^{\mathrm{ch}}$ are the molar flow rate and standard chemical exergy of component $q$, respectively. The standard chemical exergy values used in this work are listed in Table~\ref{tab:reference_chemical_exergy}.

Electricity is treated as pure exergy:
\begin{equation}
	\dot{E}_{x,\mathrm{el}}
	=
	P_{\mathrm{el}}.
	\label{eq:appendix_electrical_exergy}
\end{equation}

For heat transferred across a system boundary at temperature $T_{\mathrm{b}}$, the corresponding exergy rate is
\begin{equation}
	\dot{E}_{x,Q}
	=
	\left(
	1-\frac{T_0}{T_{\mathrm{b}}}
	\right)
	\dot{Q}.
	\label{eq:appendix_heat_exergy}
\end{equation}

Based on these definitions, the exergy balance of a control volume under quasi-steady operation follows
\begin{equation}
	\sum\nolimits_{\mathrm{in}}\dot{E}_{x,m}
	+
	\sum\nolimits_{\mathrm{in}}\dot{E}_{x,Q}
	+
	\dot{E}_{x,\mathrm{el}}^{\mathrm{in}}
	=
	\sum\nolimits_{\mathrm{out}}\dot{E}_{x,m}
	+
	\sum\nolimits_{\mathrm{out}}\dot{E}_{x,Q}
	+
	\dot{E}_{x,\mathrm{el}}^{\mathrm{out}}
	+
	\dot{E}_{x,\mathrm{d}},
	\label{eq:appendix_exergy_balance}
\end{equation}
where $\dot{E}{x,\mathrm{d}}$ is the exergy destruction caused by internal irreversibility. For heat dissipated to the surroundings at boundary temperature $T_{\mathrm{b}}$, the associated exergy loss is
\begin{equation}
	\dot{E}_{x,\mathrm{loss}}
	=
	\left(
	1-{T_0}/{T_{\mathrm{b}}}
	\right)
	\dot{Q}_{\mathrm{loss}}.
	\label{eq:appendix_exergy_loss}
\end{equation}


\section{Reference Exergy Data and Simulation Parameters of the Hydrogen Plant}
\label{sec:appendix_parameters}

\setcounter{figure}{0}
\renewcommand{\thefigure}{C\arabic{figure}}
\setcounter{table}{0}
\renewcommand{\thetable}{C\arabic{table}}

Tables~\ref{tab:awe_parameters}--\ref{tab:plant_parameters} summarize the AWE, reference chemical exergy, PTUS, control, performance-assessment, and plant-building parameters used in the simulations. Unless otherwise cited, the values are based on engineering design and operational data from the $80~\mathrm{MW}$ hydrogen plant in Northern China.

\begin{table}[H]\scriptsize\centering
	\renewcommand{\arraystretch}{1.}
	\setlength{\extrarowheight}{0pt}
	\renewcommand{\tabularxcolumn}[1]{m{#1}}
	\caption{Parameters of the AWE system.}\vspace{6pt}
	\label{tab:awe_parameters}
	\begin{tabularx}{\textwidth}{@{}>{\centering\arraybackslash}m{0.46\textwidth}>{\centering\arraybackslash}X@{}}
		\hline\hline
		Parameter & Value \\
		\hline
		Stacks per group and number of groups, $N$, $N_{\mathrm{g}}$ & $4$; $4$ \\
		Number of cells per stack, $N^{\mathrm{cell}}$ & $312$ \\
		{Reversible, thermoneutral, and maximum cell voltages, $U^{\mathrm{rev}}$, $U^{\mathrm{th}}$, $U^{\max}$} & {$1.229~\mathrm{V}$; $1.48~\mathrm{V}$; $2.1~\mathrm{V}$} \\
		Rated current, $I_{\mathrm{rate}}$ & $7{,}800~\mathrm{A}$ \\
		Maximum current ramp rate, $\Delta I_{i,j}^{\max}$ & $62.4~\mathrm{A/s}$ \\
		System pressure, $p$ & $1.8\times10^{6}~\mathrm{Pa}$ \\
		Electrochemical parameters, $r_1$, $r_2$, $r_3$, $s$ & $8.175\times10^{-6}~\Omega$; $2.136\times10^{-7}~\Omega/\mathrm{K}$; $-8.656\times10^{-12}~\Omega/\mathrm{Pa}$; $7.024\times10^{-2}$ \\
		Electrochemical parameters, $t_1$, $t_2$, $t_3$ & $-1.756\times10^{-1}~\Omega$; $79.25~\Omega\cdot\mathrm{K}$; $32.56~\Omega\cdot\mathrm{K}^{2}$ \\
		Faradaic efficiency parameters, $f_1$, $f_2$ & $50+2.5T_{i,j,\mathrm{s,out}}~(\mathrm{A}^{2})$; $0.92-6.25\times10^{-6}T_{i,j,\mathrm{s,out}}$ \\
		{Faraday constant, $F$} & {$96{,}485~\mathrm{C/mol}$} \\
		{Molar mass of hydrogen, $M_{\mathrm{H}_2}$} & {$2.016~\mathrm{g/mol}$} \\
		{Higher heating value of hydrogen, $\mathit{HHV}_{\mathrm{H}_2}$} & {$285.83~\mathrm{kJ/mol}$} \\
		{Molar volume at normal conditions, $V_{\mathrm{m,N}}$} & {$2.2414\times10^{-2}~\mathrm{m}^{3}/\mathrm{mol}$} \\
		Linearized hydrogen production parameters, $A_1$, $A_2$ & $205.31~\mathrm{Nm}^{3}/\mathrm{MWh}$; $17.85~\mathrm{Nm}^{3}/\mathrm{h}$ \\
		Standby stack power, $P_{\mathrm{SB}}$ & $0.02~\mathrm{MW}$ \\
		Minimum and maximum stack powers, $P_{\min}$, $P_{\max}$ & $1.0~\mathrm{MW}$; $5~\mathrm{MW}$ \\
		BHE and CWHE heat transfer coefficients, $k_{\mathrm{h}}$, $k_{\mathrm{c}}$ & $816.67~\mathrm{W}/(\mathrm{m}^{2}\cdot\mathrm{K})$; $980~\mathrm{W}/(\mathrm{m}^{2}\cdot\mathrm{K})$ \\
		BHE and CWHE heat exchange areas, $A_{\mathrm{h}}$, $A_{\mathrm{c}}$ & $240~\mathrm{m}^{2}$; $240~\mathrm{m}^{2}$ \\
		BHE and CWHE thermal resistances, $R_{\mathrm{BHE}}$, $R_{\mathrm{CWHE}}$ & $0.0033~\mathrm{K/W}$; $0.0033~\mathrm{K/W}$ \\
		Rated and minimum lye flow rates, $v_{\mathrm{lye,rate}}$, $\underline{v}_{\mathrm{lye}}$ & $0.0288~\mathrm{m}^{3}/\mathrm{s}$; $0.001~\mathrm{m}^{3}/\mathrm{s}$ \\
		Rated cooling water flow rate, $v_{\mathrm{cool,rate}}$ & $0.032~\mathrm{m}^{3}/\mathrm{s}$ \\
		Cooling water inlet temperature, $T_{j,\mathrm{c,in}}$ & $288~\mathrm{K}$ \\
		Lye density and specific heat capacity, $\rho_{\mathrm{lye}}$, $c_{\mathrm{lye}}$ & $1{,}250~\mathrm{kg}/\mathrm{m}^{3}$; $3{,}300~\mathrm{J}/(\mathrm{kg}\cdot\mathrm{K})$ \\
		Water densities, $\rho_{\mathrm{h}}$, $\rho_{\mathrm{c}}$ & $1{,}000~\mathrm{kg}/\mathrm{m}^{3}$; $1{,}000~\mathrm{kg}/\mathrm{m}^{3}$ \\
		Water specific heat capacities, $c_{\mathrm{h}}$, $c_{\mathrm{c}}$ & $4{,}180~\mathrm{J}/(\mathrm{kg}\cdot\mathrm{K})$; $4{,}180~\mathrm{J}/(\mathrm{kg}\cdot\mathrm{K})$ \\
		Air density and specific heat capacity, $\rho_{\mathrm{air}}$, $c_{\mathrm{air}}$ & $1.204~\mathrm{kg}/\mathrm{m}^{3}$; $1{,}005~\mathrm{J}/(\mathrm{kg}\cdot\mathrm{K})$ \\
		Group inlet temperature setpoint, $T_{\mathrm{in,set}}$ & $345~\mathrm{K}$ \\
		Stack diameter, $\varphi_{\mathrm{s}}$ & $1.90~\mathrm{m}$ \\
		Stack and SEP heat dissipation areas, $A_{\mathrm{s,diss}}$, $A_{\mathrm{sep,diss}}$ & $50.44~\mathrm{m}^{2}$; $35.0~\mathrm{m}^{2}$ \\
		Stack and SEP emissivities, $\varepsilon_{\mathrm{s}}$, $\varepsilon_{\mathrm{sep}}$ & $0.25$; $0.60$ \\
		SEP heat transfer coefficient, $h_{\mathrm{sep}}$ & $4.0~\mathrm{W}/(\mathrm{m}^{2}\cdot\mathrm{K})$ \\
		Stefan--Boltzmann constant, $\sigma$ & $5.670\times10^{-8}~\mathrm{W}/(\mathrm{m}^{2}\cdot\mathrm{K}^{4})$ \\
		Stack and SEP thermal capacities, $C_{i,j,\mathrm{s}}$, $C_{j,\mathrm{sep}}$ & $4.14\times10^{7}~\mathrm{J/K}$; $1.259\times10^{8}~\mathrm{J/K}$ \\
		BHE-side thermal capacities, $C_{j,\mathrm{he}}$, $C_{j,\mathrm{h}}$ & $1.9488\times10^{8}~\mathrm{J/K}$; $4.10\times10^{7}~\mathrm{J/K}$ \\
		CWHE-side thermal capacities, $C_{j,\mathrm{ce}}$, $C_{j,\mathrm{c}}$ & $1.9488\times10^{8}~\mathrm{J/K}$; $8.0632\times10^{7}~\mathrm{J/K}$ \\
		\hline\hline
	\end{tabularx}
\end{table}

\begin{table}[H]\scriptsize\centering
    \renewcommand{\arraystretch}{1.}
    \caption{Reference chemical exergy.}\vspace{6pt}
    \label{tab:reference_chemical_exergy}
    \begin{tabular}{ccc}
        \hline\hline
        Reference quantity & Phase & Value ($\mathrm{kJ/mol}$) \\
        \hline
        Chemical exergy of $\mathrm{H}_2(\mathrm{g})$ & Gas & $236.09$ \\
        Chemical exergy of $\mathrm{O}_2(\mathrm{g})$ & Gas & $3.97$ \\
        Chemical exergy of $\mathrm{H}_2\mathrm{O}(\mathrm{l})$ & Liquid & $0.90$ \\
        Chemical exergy of $\mathrm{H}_2\mathrm{O}(\mathrm{g})$ & Gas & $9.50$ \\
        Chemical exergy of $\mathrm{KOH}(\mathrm{s})$ & Solid & $79.40$ \\
        \hline\hline
    \end{tabular}
\end{table}

\begin{table}[H]\scriptsize\centering
    \renewcommand{\arraystretch}{1}
    \caption{Parameters of the PTUS.}\vspace{6pt}
    \label{tab:ptus_parameters}
    \begin{tabular}{ccccc}
        \hline\hline
        Branch & $v_{jk,\max}$ ($\mathrm{m}^{3}/\mathrm{s}$) & $d_{jk}$ ($\mathrm{mm}$) & $l_{jk}$ ($\mathrm{m}$) & $\psi_{jk}^{\mathrm{HL}}$ \\
        \hline
        1--2 & $0.0572$ & $200$ & $7$  & $0.9999$ \\
        2--3 & $0.0516$ & $180$ & $7$  & $0.9998$ \\
        3--4 & $0.0460$ & $160$ & $7$  & $0.9998$ \\
        4--5 & $0.0404$ & $140$ & $7$  & $0.9998$ \\
        5--6 & $0.0414$ & $120$ & $13$ & $0.9998$ \\
        6--7 & $0.0003$ & $20$  & $5$  & $0.9998$ \\
        \hline\hline
    \end{tabular}
\end{table}

\begin{table}[H]\scriptsize\centering
	\renewcommand{\arraystretch}{1.}
	\caption{Parameters of the lower-layer thermal management strategy and temperature controllers.}\vspace{6pt}
	\label{tab:thermal_management_parameters}
	\begin{tabular}{cc}
		\hline\hline
		Parameter and symbol & Value \\
		\hline
		Thermal management advance time, $t_{\mathrm{a}}$ & $9000~\mathrm{s}$ \\
		Preheating time, $t_{\mathrm{preheat}}$ & $2400~\mathrm{s}$ \\
		Thermal standby and startup target temperature, $T_{\mathrm{set}}$ & $358~\mathrm{K}$ \\
		Temperature deadband, $\Delta T_{\mathrm{d}}$ & $2~\mathrm{K}$ \\
		Low lye circulation flow rate before preheating, $v_{\mathrm{low}}$ & $0.25v_{\mathrm{lye,rate}}$ ($0.0072~\mathrm{m}^{3}/\mathrm{s}$) \\
		Lye flow rates for preheating and thermal standby, $v_{\mathrm{pre}}$, $v_{\mathrm{SB}}$ & $v_{\mathrm{lye,rate}}$ ($0.0288~\mathrm{m}^{3}/\mathrm{s}$) \\
		Rated recovery-branch flow rate, $v_{\mathrm{rec,rate}}$ & $0.0284~\mathrm{m}^{3}/\mathrm{s}$ \\
		Minimum cooling water flow rate, $\underline{v}_{\mathrm{cool}}$ & $0.05v_{\mathrm{cool,rate}}$ ($0.0016~\mathrm{m}^{3}/\mathrm{s}$) \\
		High-load current threshold, $I_{\mathrm{HL}}$ & $0.9I_{\mathrm{rate}}$ ($7{,}020~\mathrm{A}$) \\
		Group-level proportional gain, $K_{\mathrm{p,g}}$ & $0.06$ \\
		Group-level integral gain, $K_{\mathrm{i,g}}$ & $0.0667$ \\
		Group-level recovery proportional gain, $K_{\mathrm{p,rec}}$ & $0.30$ \\
		Group-level recovery integral gain, $K_{\mathrm{i,rec}}$ & $0.02$ \\
		Group-level cooling-water proportional gain, $K_{\mathrm{p,c}}$ & $0.12$ \\
		Group-level cooling-water integral gain, $K_{\mathrm{i,c}}$ & $0.90$ \\
		Stack-level proportional gain, $K_{\mathrm{p,s}}$ & $0.084$ \\
		Stack-level integral gain, $K_{\mathrm{i,s}}$ & $0.00297$ \\
		\hline\hline
	\end{tabular}
\end{table}

\begin{table}[H]\scriptsize\centering
	\renewcommand{\arraystretch}{1}
	\caption{Parameters used for performance assessment.}\vspace{6pt}
	\label{tab:awe_economic_degradation_parameters}
	\begin{tabular}{cc}
		\hline\hline
		Parameter & Value \\
		\hline
		Cold/hot-startup threshold, $T_{\mathrm{start}}$ & $333~\mathrm{K}$ \\
		Voltage degradation per cold startup, $\Delta U_{\mathrm{c}}$ & $8\times10^{-5}~\mathrm{V}$ \\
		Voltage degradation per hot startup, $\Delta U_{\mathrm{h}}$ & $8\times10^{-6}~\mathrm{V}$ \\
		Initial hydrogen conversion rate, $\gamma_{i,j,1}$ & $0.02~\mathrm{kg/kWh}$ \\
		Fatigue model parameters, $\alpha$, $\beta$ & $9000$; $1.65$ \\
		Unit CAPEX of a single stack & $8{,}300~\mathrm{CNY/kW}$ \\
		Thermal coupling topology CAPEX, $\mathrm{CAPEX}_{\mathrm{TC}}$ & $0.120\,\mathrm{CAPEX}_{i,j}$ \\
		Thermal coupling topology annual O\&M, $\mathrm{OPEX}_{\mathrm{TC},y}$ & $0.0025\,\mathrm{CAPEX}_{i,j}/\mathrm{year}$ \\
		Purchased renewable electricity price, $C^{\mathrm{W}}$ & $218.08~\mathrm{CNY/MWh}$ \\
		Hydrogen selling price, $C^{\mathrm{H}_2}$ & $2.620~\mathrm{CNY/Nm}^{3}$ \\
		Startup cost, $C^{\mathrm{SU}}$ & $802.18~\mathrm{CNY/startup}$ \\
		Shutdown cost, $C^{\mathrm{SD}}$ & $115.59~\mathrm{CNY/shutdown}$ \\
		Project lifetime, $Y$ & $15~\mathrm{years}$ \\
		Discount rate, $r$ & $5\%$ \\
		\hline\hline
	\end{tabular}
\end{table}

\begin{table}[H]\scriptsize\centering
	\renewcommand{\arraystretch}{1}
	\caption{{Parameters of the plant building and PTUS auxiliary equipment.}}\vspace{6pt}
	\label{tab:plant_parameters}
	\begin{tabular}{cc}
		\hline\hline
		Parameter & Value \\
		\hline
		Plant building volume, $V_{\mathrm{plant}}$ & $197{,}600~\mathrm{m}^{3}$ \\
		Heat transfer area of the building envelope, $A_{\mathrm{b}}$ & $15{,}200~\mathrm{m}^{2}$ \\
		Average heat transfer coefficient of the building envelope, $U_{\mathrm{b}}$ & $1.5~\mathrm{W}/(\mathrm{m}^{2}\cdot\mathrm{K})$ \\
		Equivalent thermal capacity of the indoor air and effective building thermal mass, $C_{\mathrm{plant}}$ & $1.04\times10^{9}~\mathrm{J/K}$ \\
		Radiator heat transfer area, $A_{\mathrm{rad}}$ & $16~\mathrm{m}^{2}$ \\
		Radiator thermal capacity, $C_{\mathrm{rad}}$ & $4.32\times10^{6}~\mathrm{J/K}$ \\
		Radiator water thermal capacity, $C_{\mathrm{w}}$ & $4.18\times10^{7}~\mathrm{J/K}$ \\
		Heat transfer area between water and radiator wall, $A_{\mathrm{water}}$ & $1.5~\mathrm{m}^{2}$ \\
		Convective heat transfer coefficient from water to the radiator, $U_{\mathrm{water}}$ & $500~\mathrm{W}/(\mathrm{m}^{2}\cdot\mathrm{K})$ \\
		Radiator emissivity, $\varepsilon_{\mathrm{rad}}$ & $0.90$ \\
		Number of radiators, $N_{\mathrm{rad}}$ & $18$ \\
		Maximum ventilation heat transfer coefficient, $U_{\mathrm{vent,max}}$ & $33{,}000~\mathrm{W/K}$ \\
		{Equivalent absolute PTUS pipe roughness, $K$} & {$5.0\times10^{-4}~\mathrm{m}$} \\
		{PTUS circulation pump head, $H_{\mathrm{p}}$} & {$20~\mathrm{m}$} \\
		{PTUS circulation pump efficiency, $\eta_{\mathrm{pump}}$} & {$0.75$} \\
		{Electric boiler efficiency, $\eta_{\mathrm{boiler}}$} & {$0.98$} \\
		{PTUS heat-source supply water temperature at Node~1, $T_{1}^{\mathrm{HS,s}}$} & {$363~\mathrm{K}$} \\
		\hline\hline
	\end{tabular}
\end{table}

\section*{Acknowledgments}

The authors gratefully acknowledge the financial support from the National Natural Science Foundation of China (52377116 and 52577129).

\section*{Declaration of Interest}

None.

\section*{Data Availability}

Data will be made available on request.


\end{document}